\documentclass[11pt]{article}
\usepackage{amsmath, amsfonts, amssymb}
\usepackage{graphicx, amsthm,color}

\providecommand{\U}[1]{\protect\rule{.1in}{.1in}}

\newtheorem{theorem}{Theorem}[section]
\newtheorem{lemma}[theorem]{Lemma}

\newtheorem{corollary}[theorem]{Corollary}
\newtheorem{proposition}[theorem]{Proposition}
\newtheorem{definition}[theorem]{Definition}
\theoremstyle{definition}

\newtheorem{example}[theorem]{Example}
\theoremstyle{remark}
\newtheorem{remark}[theorem]{Remark}
\numberwithin{equation}{section}

\newcommand{\R}{\mathbb{R}}

\newcommand{\N}{\mathbb{N}}

\providecommand{\dom}{\mathop{\rm dom}\nolimits}
\providecommand{\argmin}{\mathop{\rm argmin}}
\providecommand{\argmax}{\mathop{\rm argmax}}

\usepackage{tikz}
\usepackage{tikz-3dplot}
\usetikzlibrary{calc,patterns,angles,quotes}
\usetikzlibrary{babel}
\usepackage{pgfplots}
\usepgfplotslibrary{patchplots}
\usetikzlibrary{backgrounds, intersections}
\usepgfplotslibrary{fillbetween}
\usetikzlibrary{arrows,automata}
\usepackage{subcaption}

\pgfplotsset{compat=1.18}
\usepgfplotslibrary{colorbrewer}
\usepackage{filecontents}

\definecolor{forestgreenweb}{rgb}{0.13, 0.55, 0.13}

\begin{document}

\begin{center}
{\LARGE Metric compatibility and determination in complete metric spaces}

\vspace{1cm}

{\Large \textsc{Aris Daniilidis, Tri Minh Le, David Salas}}
\end{center}

\bigskip

\noindent\textbf{Abstract.} It was established in \cite{DS2022} that Lipschitz
inf-compact functions are uniquely determined by their local slope and
critical values. Compactness played a paramount role in this result, ensuring
in particular the existence of critical points. We hereby emancipate from this
restriction and establish a determination result for merely bounded from below
functions, by adding an assumption controlling the asymptotic behavior. This
assumption is trivially fulfilled if $f$ is inf-compact. In addition, our
result is not only valid for the (De Giorgi) local slope, but also for the
main paradigms of average descent operators as well as for the global slope,
case in which the asymptotic assumption becomes superfluous. Therefore, the
present work extends simultaneously the metric determination results of
\cite{DS2022} and \cite{TZ2022}. \bigskip


\vspace{0.6cm}

\textit{Secondary} 03F15 ; 30L15

\tableofcontents


\section{Introduction}


\noindent In \cite[Theorem~2.4]{DS2022} the authors have shown that in every
metric space, the local slope operator contains sufficient information to
determine any continuous inf-compact function with finite slope. Indeed,
knowledge of the critical values (values of the function on the set of points
where the slope is zero) and knowledge of the slope at every point determine
uniquely the function. We hereafter refer to this result as
\emph{determination result}. The proof makes use of transfinite induction and
is based on a cardinality obstruction. Pertinence of the assumptions was also
thereby discussed. \smallskip\newline In the follow-up work \cite{DMS2022} the
authors adopted a much more general framework: they introduced an abstract
notion of descent modulus, based on three axioms (see \cite[Definition~3.1]
{DMS2022} or properties $(\mathcal{D}_{1})$--$(\mathcal{D}_{3})$ of
forthcoming Definition~\ref{def:T}) and showed that the result of
\cite{DS2022} can be emancipated from the metric structure and fit to a mere
topological setting, provided a reasonable notion of \textit{steepest descent}
(or other meaningful notion of descent, like \textit{average descent}) is
coined. Therefore, instead of considering metric spaces, we can work on
probability spaces or Markov chains. However, similarly to~\cite{DS2022}, an
underlying compactness assumption was still required in~\cite{DMS2022}: the
functions for which the result applies should (be continuous and) have compact
sublevel sets. This was indeed paramount for the proof of the main result of
both works. \smallskip\newline The aim of the current work is to eliminate the
compactness assumption and use instead completeness together with a control on
asymptotic behaviour. This renounces full generality, restricting naturally to
the framework of (complete) metric spaces.\smallskip\newline Very recently, in
the same setting of complete metric spaces, Thibault and Zagrodny in
\cite{TZ2022} were able to obtain a determination result for the
\textit{global slope} (we recall this definition in~\eqref{eq:globalSlope}).
The proof is highly technical and uses the notion of countably orderable families
previously introduced in~\cite{GZ1992}. For a general function, the global
slope is a very restrictive notion (controlling also the asymptotic behavior),
but for the class of convex functions it coincides with the local slope and
the authors were able to obtain the following powerful convex determination result:

\begin{itemize}
\item (convex determination) \textit{Two convex continuous and
bounded from below functions \newline with the same slope can only differ by a
constant.}
\end{itemize}

The above result was initially established in Hilbert spaces,
see~\cite{BCD2018} (smooth case) and~\cite{PSV2021} (nonsmooth case). It can
also be obtained as a corollary of a more general sensitivity result, derived
in~\cite{DD2023}, which states, roughly speaking, that the slope deviation
between two convex functions controls the deviation between the functions
themselves. A similar determination result was obtained using proximal operators \cite{V2021}. All these proofs rely heavily on (sub)gradient descent systems,
making crucial use of the Hilbertian structure. However, this drawback no
longer appears in~\cite{TZ2022}, where the authors, working directly in metric
spaces with the global slope, were able to establish the validity of the above
convex determination result in Banach spaces.\smallskip\newline Coming back to
the present work, we enhance the technique developed in~\cite{DS2022} to
obtain a general determination result in the setting of complete metric
spaces. Comparing with~\cite{TZ2022}, the result not only applies for the
global slope (where the interest is essentially limited to the convex
determination in a Banach space), but also for the local slope (the definition
is recalled in~\eqref{eq:metricSlope}) as well as for the main paradigms of
average descent operator discussed in Section~\ref{ssec:2.3}. As a
consequence, the result applies to a large class of functions (for instance,
Lipschitz functions in complete metric spaces). This already hints potential
applications in Eikonal equations, or more generally, in Hamilton-Jacobi
equations whose viscosity solutions admit an alternative description via
slopes (see \cite{GS2015}, \cite{GHN2015}, \cite{LSZ2021} \emph{e.g.}).
A~further extension is made by formulating the result in terms of an abstract
descent modulus in the spirit of~\cite{DMS2022}, but with an extra property
(metric compatibility) to reckon with the given metric (see
Definition~\ref{def:Tmetric}). From a practical viewpoint, in a given metric
space all reasonable descent moduli are metrically compatible (see also
discussion in Section~\ref{ssec:2.4}).


\subsection{Organization of the manuscript}

The manuscript is organized as follows: in Subsection~\ref{ssec:1.2} we fix
notation and terminology, while in Section~\ref{sec:pre} we revisit
from~\cite{DMS2022} the definition of an abstract descent modulus and readjust
it (\emph{c.f.} Definition~\ref{def:T}) to encompass extended real-valued
functions in a way that the determination result still holds for inf-compact
functions which are continuous in their domain.\smallskip\newline
Subsection~\ref{ssec:2.1} resumes the State--of--the--art in this (slightly)
more general setting, with the extra benefit that the proofs are now
significantly simplified. This is possible because Definition~\ref{def:T} is
defined in a compatible way with respect to function truncation, see proofs of
Lemma~\ref{lemma_sal1} and Theorem~\ref{thm:DetExtended-compact}. We then
obtain Corollary~\ref{cor:inf-compact} which readily
extends~\cite[Proposition~2.2]{DS2022} and~\cite[Theorem~2.4]{DS2022}
.\smallskip\newline In Subsection~\ref{ssec:2.2} we establish an easy
noncompact determination result for the case of smooth functions in a Banach
space for the natural descent modulus $T[f]=\lVert\nabla f\lVert$. The result illustrates
perfectly the need of controlling the asymptotic behaviour and at the same
time hints towards the right definition of \textit{asymptotically critical
sequence} (see~Definition~\ref{def:AC-1}).\smallskip\newline In
Subsection~\ref{ssec:2.3} we present the main paradigms of descent in a metric
space which are covered by our main result: the (De Giorgi's) local slope, the
global slope, the average descent and the diffusion descent. These paradigms
are recalled in Subsection~\ref{ssec:2.4} and treated in a uniform manner by
means of the definition-scheme of a \textit{metrically compatible descent
modulus} (Definition~\ref{def:Tmetric}).\smallskip\newline The main result is
presented in Section~\ref{sec:3}. Controlling the critical values, the
asymptotic behavior and the abstract descent at each point leads to
Theorem~\ref{thm.strict_compa} (comparison lemma) and
Theorem~\ref{thm:DetExtended-metric} (determination result). We recover the
determination result of~\cite{TZ2022} as a corollary, by applying our result
for the global slope, which is a particular case of an abstract descent, since
in this case every asymptotically critical sequence is infimizing for the function.


\subsection{Notation and terminology.}

\label{ssec:1.2} Throughout this work $X$ is a complete metric space, which
will be eventually upgraded to a Banach space in Subsection~\ref{ssec:2.2}.
Given any function $f:X\rightarrow\mathbb{R}$ and $r\in\mathbb{R}$ we define
by
\[
\lbrack f\leq r]:=\{x\in X:\,f(x)\leq r\},
\]
the $r$-sublevel set. The strict sublevel set $[f<r]$ is defined analogously.
We denote the (effective) domain of $f$ by
\[
\dom f:=\{x\in X:\,f(x)<+\infty\}.
\]
\smallskip\newline A~function $f$ is lower semicontinuous (in short, lsc) if
for every $r\in\mathbb{R}$ the sublevel set $[f\leq r]$ is closed. We say that
$f$ is proper if it has at least one nonempty sublevel set, or equivalently,
if $\dom f\neq\emptyset$. Further, a function $f$ is called \emph{inf-compact}
if the sublevel sets $[f\leq r]$ are compact for all $r<\sup f$. Notice that
every lower semicontinuous inf-compact function attains its global
minimum.\smallskip\newline We further denote by $(\mathbb{R}\cup
\{+\infty\})^{X}$ the set of extended real-valued functions on $X$ and by
$\mathcal{C}(X)$ the set of continuous real-valued functions on $X.$ We also
denote by
\begin{align*}
\mathcal{LSC}(X)  &  :=\left\{  f:X\rightarrow\mathbb{R}\cup\{+\infty
\}\ :\ f\text{ proper, lsc}\right\}  ;\\
\overline{\mathcal{C}}(X)  &  :=\left\{  f:X\rightarrow\mathbb{R}\cup
\{+\infty\}\ :\ f\text{ proper, lsc and }f|_{\dom f}\text{ continuous}
\right\}  .
\end{align*}
A subset $\mathcal{F}$ of $(\mathbb{R}\cup\{+\infty\})^{X}$ is called a
\emph{cone}, if for any $f\in\mathcal{F}$ and $r\geq0$ we have $rf\in
\mathcal{F}$ (with the convention $0\cdot(+\infty)=0$). A cone $\mathcal{F}$
which is closed under translation (that is, $f+c\in\mathcal{F}$ for all
$f\in\mathcal{F}$ and $c\in\mathbb{R}$) will be called \emph{translation
cone}. In what follows, $\mathcal{F}$ will always be a translation cone of
proper functions.\medskip\newline Notice that $\mathcal{C}(X)$ is a vector
subspace of $(\mathbb{R}\cup\{+\infty\})^{X}$ while $\mathcal{LSC}(X)$,
$\overline{\mathcal{C}}(X)$ are translation cones. Notice further that
$f\in\mathcal{C}(X)$ if and only if $f\in\overline{\mathcal{C}}(X)$ and
$\dom
f=X.$\medskip\newline For any $a\in\mathbb{R}$ we set $a^{+}:=\max\{a,0\}$
(the positive part of the number $a$).\smallskip\newline For a function
$f:X\rightarrow\mathbb{R}\cup\{+\infty\}$ we define the operator $\Delta
^{+}f:\dom f\times X\rightarrow\mathbb{R}$ by
\begin{equation}
\Delta^{+}f(x,y)=\left\{
\begin{array}
[c]{cl}
\frac{\{f(x)-f(y)\}^{+}}{d(x,y)}\,,\smallskip & \text{ if }x\neq y,\smallskip\\
0\,, & \text{ if }x=y.
\end{array}
\right.  \label{eq:DeltaPlus}
\end{equation}


\section{Descent moduli: state-of-the-art and extended definitions}

\label{sec:pre}

Following the spirit of \cite{DMS2022}, we call \textit{descent modulus} on a
topological space $X$ any operator $T:\mathcal{F}\rightarrow\lbrack
0,+\infty]^{X}$ satisfying three natural properties (see $(\mathcal{D}_{1}
)$--$(\mathcal{D}_{3})$ in Definition~\ref{def:T} below). The quantity
$T[f](x)\in\lbrack0,+\infty]$ is conceived as an abstract measurement of
descent for the function $f$ at the point $x$. If $T[f](x)=0,$ then the point
$x$ is called $T$-\textit{critical} (or simply critical). Therefore, the set
of $T$-critical points of $f$ coincides with the zeros of the function $T[f]$
and is denoted by
\begin{equation}
\mathcal{Z}_{T}(f):=\{x\in X:\,T[f](x)=0\}. \label{eq:crit}
\end{equation}

A formal definition for proper extended real-valued functions follows:

\begin{definition}
[Descent modulus]\label{def:T}Let $\mathcal{F}\subset(\mathbb{R}\cup
\{+\infty\})^{X}$ be a translation cone.

An operator $T:\mathcal{F}\rightarrow\lbrack0,+\infty]^{X}$ is called descent
modulus on $\mathcal{F}$ if
\begin{equation}
\dom T[f]\subset\dom f,\text{\quad for every }f\in\mathcal{F} \tag{$\mathcal{D}_0$}
\end{equation}
and the following three conditions hold for every $f,g\in\mathcal{F}$ and
$x\in X:$
\end{definition}

\begin{itemize}
\item[$(\mathcal{D}_{1})$] $x\in\argmin f\;\implies\;x\in\mathcal{Z}_{T}(f).$

\item[$(\mathcal{D}_{2})$] $T[f](x)<T[g](x)\;\implies\;\exists z\in\dom
g:\,\{f(x)-f(z)\}^{+}<\{g(x)-g(z)\}^{+}.$

\item[$(\mathcal{D}_{3})$] If $0<T[f](x)<+\infty$ and $r>1,$ then
$T[f](x)<T[rf](x).$
\end{itemize}

Let us have a brief discussion on the properties defining the descent modulus:
Property $(\mathcal{D}_{1})$ guarantees preservation of global minima, since
$\mathcal{Z}_{T}(f)=\argmin T[f].$ Property $(\mathcal{D}_{2})$ can be seen as
a monotonicity property on the sublevel set: indeed, if for every $z\in\lbrack
g\leq g(x)]$ one has $f(x)-f(z)\geq g(x)-g(z)\,$(that is, if $f$ has more
descent than $g$ in all descent directions of $g$) then one should necessarily
have $T[f](x)\geq T[g](x).$ Therefore, $(\mathcal{D}_{2})$ can be restated as
follows:
\begin{equation}
\left.
\begin{array}
[c]{c}
\left\{  g(x)-g(z)\right\}  ^{+}\leq\left\{  f(x)-f(z)\right\}  ^{+}
\smallskip\\
\text{for all }z\in\dom g
\end{array}
\right\}  \,\Longrightarrow\,T[g](x)\leq T[f](x). \tag{$\widetilde{\mathcal
{D}_2}$}
\end{equation}
Finally $(\mathcal{D}_{3})$ is a scalar monotonicity property, ensuring that
if a function $f$ has a nonzero finite descent at $x,$ then the function
$(1+\varepsilon)f$ has an amplified descent at the same point for any
$\varepsilon>0$.

\begin{remark}
\textrm{(i).} Definition~\ref{def:T} applies to extended real-valued functions
and $(\mathcal{D}_{0})$ imposes an infinite descent to all points for which
$f(x)=+\infty.$ If $\mathcal{F}=\mathcal{C}(X),$ then $(\mathcal{D}_{0})$
holds trivially and the above definition of descent modulus coincides with
\cite[Definition 3.1]{DMS2022}.\smallskip\newline\textrm{(ii).} A
straightforward consequence of the above definition is that a descent modulus
can be defined only for proper functions. (To see this, given $g\in
\mathcal{F}$, consider the function $f\equiv\boldsymbol{0}$ in $\mathcal{F}$
and apply $(\mathcal{D}_{2})$).
\end{remark}

We define the domain $\dom T\subset\mathcal{F}$ of a descent modulus $T$ as
follows:
\begin{equation}
\dom T:=\{f\in\mathcal{F}:\,\dom T[f]=\dom f\}, \label{domT}
\end{equation}
that is, $f\in\dom T$ if and only if it has a finite slope at every point in
which it has a finite value.

If $X$ is a metric space, then $\dom T$ contains the class of Lipschitz
continuous functions for every reasonable descent modulus. (The reader can
easily verify that this is the case for the main instances of descent moduli
of this work: \textit{c.f.} Example~\ref{ejem32} and Example~\ref{ejem26}.)


\subsection{Determination in compact spaces}

\label{ssec:2.1} The determination result established in~\cite{DMS2022}
requires the functions to have compact sublevel sets. The proof was based on a
transfinite induction and the conclusion was obtained by contradiction, due to
a cardinality obstruction since the induction did not allow point repetitions.
In this section, for the sake of completeness, we restate this result in a
slightly more general setting: the descent modulus is now considered on
extended real-valued (inf-compact) functions. In fact, this new framework,
contemplated by the extended Definition~\ref{def:T} allows a much simpler
proof (namely, the transfinite induction is replaced by a maximum principle),
which in the setting of~\cite{DS2022}, \cite{DMS2022} was formally impossible.
We present this proof here.

\begin{lemma}
[Strict comparison in compact spaces]\label{lemma_sal1} Let $X$ be a compact
topological space and $T$ a descent modulus on a translation cone
$\mathcal{F}$ containing $\mathcal{LSC}(X)$. Let $f\in\overline{\mathcal{C}
}(X)$, $g\in\mathcal{LSC}(X)$ and assume:

\begin{itemize}
\item[$\mathrm{(i).}$] (descent domination) $T[f](x)<T[g](x)$, for every
$x\in\dom g\setminus\mathcal{Z}_{T}(g)$ ;

\item[$\mathrm{(ii).}$] (control of criticality) $f(z)<g(z),$ for every
$z\in\mathcal{Z}_{T}(g)$ ;
\end{itemize}

Then, it holds
\[
f(x)<g(x),\qquad\forall x\in\dom g.
\]

\end{lemma}

\noindent\textbf{Proof.} Notice that $\dom g\subset\dom f$, therefore $f$ is
continuous on $\dom g$. Let us first assume that $g$ is finite, that is,
$\dom
g=X$. Then, $f-g$ is (finite and) upper semicontinuous and attains its maximum
at some point $x_{0}\in X$. It suffices to show that $x_{0}\in\mathcal{Z}
_{T}(g)$ (then~(ii) applies and $\max\,(f-g)=(f-g)(x_{0})<0$). If $x_{0}
\notin\mathcal{Z}_{T}(g)$, then, $T[g](x_{0})>T[f](x_{0})$ which yields by
hypothesis $(\mathcal{D}_{2})$ that there exists $z\in X$ such that
$\{f(x_{0})-f(z)\}^{+}<\{g(x_{0})-g(z)\}^{+}$. In particular, $(f-g)(x_{0}
)<(f-g)(z)$, which is a contradiction.\smallskip\newline Let us now consider
the case $\dom g\neq X$, that is, $g$ takes the value $+\infty$ at some point.
Let $h:X\rightarrow\mathbb{R}\cup\{+\infty\}$ given by
\[
h(x)=
\begin{cases}
(f-g)(x)\,,\quad & \text{ if }x\in\dom g\,,\\
\phantom{jo}+\infty\,, & \text{ otherwise.}
\end{cases}
\]
Fix any $a>\inf g$. Since $g$ is lsc, the set $K_{a}=[g\leq a]$ is nonempty
compact and the upper semicontinuous function $h$ attains its maximum there at
some point $x_{a}\in K_{a}$. If $x_{a}\notin\mathcal{Z}_{T}(g)$, then, as
before, there exists $z_{a}\in\dom g$ such that $\{f(x_{a})-f(z_{a}
)\}^{+}<\{g(x_{a})-g(z_{a})\}^{+}$. This yields that $z_{a}\in K_{a}$ and
$h(x_{a})<h(z_{a})$, which is a contradiction. Thus, $x_{a}\in\mathcal{Z}
_{T}(g)$ and $h$ is strictly negative in $K_{a}$. Since $\dom g=\bigcup
_{a>\inf g}[g\leq a]$, the conclusion follows.\hfill$\Box$

\bigskip

The following theorem is the direct extension of the determination theorems
of~\cite{DMS2022}, invoking Lemma~\ref{lemma_sal1} instead of~\cite[Lemma~3.3]
{DMS2022}.

\begin{theorem}
[Descent determination of extended real-valued functions in compact
spaces]\label{thm:DetExtended-compact} Let $X$ be a compact topological space
and $T$ a descent modulus on a translation cone $\mathcal{F}$ containing
$\mathcal{LSC}(X)$. Let $f\in\overline{\mathcal{C}}(X)$ and $g\in
\mathcal{LSC}(X) \cap\dom(T)$. Then,

\begin{enumerate}
\item[(a)] If $T[f](x)\leq T[g](x)$, for all $x\in X$ and $f(x)\leq g(x)$, for
all $x\in\mathcal{Z}_{T}(g)$, then $f\leq g$.

\item[(b)] If $f,g\in\overline{\mathcal{C}}(X)\cap\dom(T)$, $T[f](x)=T[g](x)$,
for all $x\in X$ and $f(x)=g(x)$ for all $x\in\mathcal{Z}_{T}(g)=\mathcal{Z}
_{T}(f)$, then $f=g$.
\end{enumerate}
\end{theorem}

\noindent\textbf{Proof.} Since statement $(b)$ is symmetric, it is sufficient
to prove $(a)$. To this end, replacing $g$ by $g_{\varepsilon}=(1+\varepsilon
)(g+\varepsilon)$, we get that $\dom
g_{\varepsilon}=\dom g$, $T[f](x)<T[g_{\varepsilon}](x)$ for every $x\in\dom
g_{\varepsilon}$, $\mathcal{Z}_{T}(g_{\varepsilon})\subset\mathcal{Z}_{T}(g)$
and $f(x)<g_{\varepsilon}(x)$ for every $x\in\mathcal{Z}_{T}(g_{\varepsilon})
$. Thus, $f<g_{\varepsilon}$ over $\dom g_{\varepsilon}$. Taking
$\varepsilon\rightarrow0$, we obtain $f\leq g$ on $\dom g$, which readily
yields $f\leq g$ on the whole space $X$. \hfill$\Box$

\bigskip

The comparison principles and the determination result of \cite{DMS2022} for
inf-compact functions can be derived from Lemma~\ref{lemma_sal1} and
Theorem~\ref{thm:DetExtended-compact}, as the following corollary shows.

\begin{corollary}
\label{cor:inf-compact}Let $X$ be a topological space (not necessarily
compact) and $T$ a descent modulus on $\mathcal{F} = \overline{\mathcal{C}
}(X)$. Let $f,g\in\mathcal{F}$ be bounded from below.

\begin{enumerate}
\item[(a)] If $g$ is inf-compact, $T[f](x)<T[g](x)$ for all $x\in
X\setminus\mathcal{Z}_{T}(g)$, and $f(x) < g(x)$ for all $x\in\mathcal{Z}
_{T}(g)$, then $f<g$.

\item[(b)] If $g$ is inf-compact, $g\in\dom(T)$, $T[f](x)\leq T[g](x)$, for
all $x\in X$ and $f(x)\leq g(x),$ for all $x\in\mathcal{Z}_{T}(g)$, then
$f\leq g$.

\item[(c)] If $f,g$ are inf-compact, $f,g\in\overline{\mathcal{C}}
(X)\cap\dom(T)$, $T[f](x)=T[g](x)$, for all $x\in X$ and $f(x)=g(x)$ for all
$x\in\mathcal{Z}_{T}(g)=\mathcal{Z}_{T}(f)$, then $f=g$.
\end{enumerate}
\end{corollary}

\noindent\textbf{Proof.} It is enough to prove $(a)$. The conclusion is
trivial if $g$ is constant, since in this case $\mathcal{Z}_{T}(g)=X$.
Therefore we may assume that $\inf g<\sup g$. Fix any $a\in(\inf g,\sup g)$
and set $K_{a}=[g\leq a].$ Then, $K_{a}$ is nonempty and compact. We set
$\mathcal{F}_{a}:=\{h\in\mathcal{F}:\,K_{a}\subset\dom h\}$ and for every
$h\in\mathcal{F}$, we define
\[
h_{a}:=h+i_{K_{a}},
\]
where $i_{K_{a}}$ denotes the indicator function of $K_{a},$ that is,
\[
i_{K_{a}}(x):=\left\{
\begin{array}
[c]{cc}
0, & x\in K_{a}\\
+\infty, & x\notin K_{a}.
\end{array}
\right.
\]
Notice that $h_{a}\in\overline{\mathcal{C}}(X)=\mathcal{F}$ and $\mathcal{F}
_{a}\subset\mathcal{F}$, so the operator $T$ is a descent modulus on
$\mathcal{F}_{a}$. We deduce easily from $(\widetilde{\mathcal{D}_{2}})$ that
for every\ $h\in\mathcal{F}_{a}$
\[
T[h_{a}](x)\leq T[h](x),\qquad\text{for all }x\in K_{a}
\]
and (since $K_{a}$ is a nontrivial sublevel set of $g$)
\[
\mathcal{Z}_{T}(g_{a})=\mathcal{Z}_{T}(g)\cap K_{a}\qquad\text{and}\qquad
T[g_{a}](x)=T[g](x),\quad\text{for all }x\in K_{a}.
\]
Then,
\[
T[f_{a}](x)\leq T[f](x)<T[g](x)=T[g_{a}](x).
\]
Moreover, for every $x\in\mathcal{Z}_{T}(g_a)$ we have $f_{a}(x)=f(x)<g(x)=g_{a}
(x)$. Applying Lemma~\ref{lemma_sal1} we deduce that $f_{a}<g_{a}$ on $K_{a}$,
and consequently, $f(x)<g(x)$, for every $x\in K_{a}$. Since $a\in(\inf g,\sup
g)$ is arbitrary, we conclude that $f<g$ on $\dom g\setminus\argmax g$.
\smallskip\newline If $\argmax g=\emptyset$ or $\sup g=+\infty$, we have
$\dom
g\setminus\argmax g=\dom g$ and the result follows. Thus, it suffices to
consider the case $\sup g<+\infty$ and $\argmax g\neq\emptyset$. Then take
$x\in\argmax g$. If $T[g](x)=0$, then $f(x)<g(x)$ by hypothesis. If not,
$T[f](x)<T[g](x)$ and property $(\mathcal{D}_{2})$ entails that there exists
$z\in\dom g$ such that
\[
f(x)-f(z)\leq\{f(x)-f(z)\}^+<\{g(x)-g(z)\}^{+}=g(x)-g(z).
\]
Then, $z\notin\argmax g$, entailing that $g(z)>f(z)$. Thus,
$f(x)<g(x)+f(z)-g(z)<g(x)$. We conclude that, regardless the value of
$T[g](x)$, we always have that $f(x)<g(x)$. Therefore $f<g$ on $\dom g$, and
the proof is complete. \hfill$\Box$


\subsection{A simple noncompact result: the smooth case}

\label{ssec:2.2} Let us consider the setting given by $X = \mathbb{R}^{d}$,
the cone $\mathcal{F} = \mathcal{C}^{1}(\mathbb{R}^{d})$ of continuously
differentiable functions, and the descent modulus given by $T[f](x) = \|\nabla
f(x)\|$.

In this setting, consider two functions $f,g\in\mathcal{F}$ bounded from below
such that
\[
\Vert\nabla f(x)\Vert\leq\lVert\nabla g(x)\lVert,\quad\forall x\in X.
\]
Following~\cite{PSV2021}\footnote{A first version of this idea was due to
J.-B. Baillon in 2018, see \cite{penot}.} we compare the functions $f$ and $g$
along the descent curves of $g$ (which is the function with dominating slope).
Given $x_{0}$, we consider the curve $\gamma:[0,+\infty)\rightarrow
\mathbb{R}^{d}$ that solves
\begin{equation}
\begin{cases}
\dot{\gamma}(t)=-\nabla g(\gamma(t)),\quad t\geq0,\\
\gamma(0)=x_{0}.
\end{cases}
\label{eq:GS}%
\end{equation}
Then, we can directly write
\begin{equation}
\begin{aligned} g(x_0) - f(x_0) &= \limsup_{t\to +\infty}\,\, (g-f)(\gamma(t)) - \int_{0}^t ((g-f)\circ \gamma)'(s) \,ds\\ &=\, \limsup_{t\to +\infty}\,\, (g-f)(\gamma(t)) - \int_{0}^t \Bigl( \langle \nabla g(\gamma(s)),\, \dot{\gamma}(t)\rangle - \langle \nabla f(\gamma(s)),\, \dot{\gamma}(t)\rangle\Bigr) \,ds \\ &= \limsup_{t\to +\infty}\, (g-f)(\gamma(t))+ \int_{0}^t \lVert\nabla g(\gamma(s))\lVert\cdot |\dot{\gamma}(s)|\, ds + \int_{0}^t \langle \nabla f(\gamma(s)),\dot{\gamma}(s)\rangle ds \\ &\geq \limsup_{t\to +\infty}\,\, (g-f)(\gamma(t)) + \int_{0}^t \Bigl( \underbrace{\lVert\nabla g(\gamma(s))\lVert - \lVert\nabla f(\gamma(s))\lVert}_{\geq 0} \Bigr) |\dot{\gamma}(s)| \,ds\\ &\geq \limsup_{t\to +\infty}\,\, (g-f)(\gamma(t)) \, = \, \lim_{t\to +\infty}\,\, g(\gamma(t)) - \liminf_{t\to +\infty}\,\, f(\gamma(t)) . \end{aligned} \label{eq:InequalitySteepestDescent}
\end{equation}
Consequently, if there exists $\bar{x}\in\omega\text{-}\lim\gamma$ ($\omega
$-limit point of $\gamma$), then $\bar{x}\in Z:=\mathcal{Z}_{T}(g)$. In such a
case, if $g\big|_{Z}\geq f\big|_{Z}$ (which is the boundary condition of
Theorem~\ref{thm:DetExtended-compact}\thinspace($a$)), we would conclude that
$g(x_{0})\geq f(x_{0})$. However, since the space is noncompact,
$\omega\text{-}\lim\gamma$ might be empty and we need to include the boundary
condition
\begin{equation}
\liminf_{t\rightarrow+\infty}\,\,f(\gamma(t))\leq\lim_{t\rightarrow+\infty
}\,\,g(\gamma(t)). \label{eq:sal}
\end{equation}
The above condition should be imposed only for steepest descent curves
$\gamma:[0,+\infty)\rightarrow\mathbb{R}^{d}$ of $g$ without $\omega$-limits
and not for those for which $\omega\text{-}\lim\gamma\neq\emptyset,$ since
this latter case is already captured by the comparison condition on the
critical set $\mathcal{Z}_{T}(g)$ ($\omega$-limits of the gradient descent
curve are automatically critical points for $g$).\smallskip\newline The
drawback of the boundary condition (\ref{eq:sal}) is that it depends on
$\nabla g$, via (\ref{eq:GS}), rather than on the descent modulus
$T[g]=\lVert\nabla g\lVert$. To overcome this difficulty and obtain a boundary
condition that is independent of $\nabla g$ we introduce the following definition:

\begin{definition}
[Asymptotically critical path, smooth version]\label{def:salas}We say that a
differentiable curve $\tilde{\gamma}:[0,+\infty)\rightarrow\mathbb{R}^{d}$ is
an asymptotically $T$-critical path for the function $g$ (with $T[g]=\lVert
\nabla g\lVert$) if
\begin{equation}
\Vert\tilde{\gamma}^{\prime}\Vert=1,\qquad\omega\text{-}\lim\tilde{\gamma
}=\emptyset\qquad\text{and }\qquad\int_{0}^{+\infty}T[g](\tilde{\gamma
}(s))\,ds<+\infty. \label{eq:jpp}
\end{equation}
\end{definition}

The above definition encompasses the following key elements for the class of
continuously differentiable bounded from below functions:

\begin{itemize}
\item[(1)] If 
$\tilde{\gamma}$ is an asymptotically critical path for $g$,
then
$\underset{s\rightarrow\infty}{\lim}T[g](\tilde{\gamma}(s))=0.$

\item[(2)] Every steepest descent curve without $\omega$-limits yields, upon
reparametrization, an asymptotically critical path.

\item[(3)] If two continuously differentiable functions have the same slope, then they have the same
critical points and the same asymptotically critical paths.
\end{itemize}

The last assertion is obvious from Definition~\ref{def:salas}, while the first
is a straightforward consequence of the integrability condition in
(\ref{eq:jpp}). Concerning the second assertion, let $\gamma:[0,+\infty
)\rightarrow\mathbb{R}^{d}$ be a descent curve of $g$ without $\omega$-limits.
Then, $\gamma$ has infinite length and so does $\tilde{\gamma}$, its arc-length
parametrization, defined by
\begin{equation}
\tilde{\gamma}(s):=\gamma(\sigma^{-1}(s))\qquad\text{where}\qquad
s=\sigma(t):= {\displaystyle\int\limits_{0}^{t}} \,\Vert\dot{\gamma}(t)||\,dt.
\label{eq:jpp1}
\end{equation}
Then, $\omega$-$\lim\tilde{\gamma}=\emptyset$ and $\Vert\tilde{\gamma}^{\prime
}(s)\Vert=1.$ Performing a change of variables we deduce:
\begin{equation}
g(\gamma(0))-\inf g\,\geq\,-\int_{0}^{+\infty}(g\circ\gamma)^{\prime
}(t)\,dt\,=\,\int_{0}^{+\infty}\lVert\nabla g(\gamma(t))\lVert\cdot\Vert
\dot{\gamma}(t)\lVert dt\,=\,\int_{0}^{+\infty}\lVert\nabla g(\tilde{\gamma
}(s))\lVert\,ds, \label{eq:jpp2}
\end{equation}
which yields (\ref{eq:jpp}).

With the above in mind, we establish the following noncompact determination result.

\begin{theorem}
\label{thm:Baillon} Let $f,g:\mathbb{R}^{d}\rightarrow\mathbb{R}$ be two
continuously differentiable functions which are bounded from below. Assume that

\begin{enumerate}
\item[(i)] $\lVert\nabla f(x)\lVert=\lVert\nabla g(x)\lVert$ for every
$x\in\mathbb{R}^{d}$;

\item[(ii)] $f(z)=g(z)$ for every $z\in\mathcal{Z}_{T}(f)=\mathcal{Z}_{T}(g)$.

\item[(iii)] $\displaystyle\liminf_{s\rightarrow+\infty}f(\tilde{\gamma
}(s))=\displaystyle\liminf_{s\rightarrow+\infty}g(\tilde{\gamma}(s))$, for
each (common) asymptotically critical path $\tilde{\gamma}$.
\end{enumerate}

Then, $f=g$.
\end{theorem}

\noindent\textbf{Proof.} Take $x_{0}\in\mathbb{R}^{d}$ and $\gamma$ be a
steepest descent curve of $g$ emanating from $x_{0}$. Then, we can divide our
analysis in two cases:\medskip

\textit{Case 1}: $\omega\text{-}\lim\gamma$ is nonempty. Choose $\bar{z}
\in\omega\text{-}\lim\gamma$. Since $\underset{t\rightarrow\infty}{\lim}
\lVert\nabla g(\gamma(t))\lVert=0,$ continuity of $\nabla g$ entails that
$\bar{z}\in\mathcal{Z}_{T}(g)$. Moreover, continuity of $g$ and $f$ yield that
$\underset{t\rightarrow\infty}{\lim}g(\gamma(t))=g(\bar{z})$ and
$\underset{t\rightarrow\infty}{\lim}f(\gamma(t))=f(\bar{z})$. Then,
\eqref{eq:InequalitySteepestDescent} yields that
\[
g(x_{0})-f(x_{0})\geq g(\bar{z})-f(\bar{z})=0.
\]

\textit{Case 2}: $\omega\text{-}\lim\gamma$ is empty.\medskip

Then, since $\gamma$ is the steepest descent curve of $g$ emanating from
$x_{0}$, we get that
\[
\int_{0}^{+\infty}\lVert\nabla g(\gamma(t))\lVert\cdot\Vert\dot{\gamma
}(t)\lVert dt=-\int_{0}^{+\infty}(g\circ\gamma)^{\prime}(t)dt=g(x_{0}
)-\lim_{t\rightarrow+\infty}g(\gamma(t))\leq g(x_{0})-\inf g<+\infty.
\]
Thus, the curve $\tilde{\gamma}\ $given in (\ref{eq:jpp1}) is an
asymptotically critical path and
consequently assumption (iii) applies. We deduce from
(\ref{eq:InequalitySteepestDescent}) that
\[
g(x_{0})-f(x_{0})\geq\lim_{t\rightarrow+\infty}g(\gamma(t))-\liminf
_{t\rightarrow+\infty}f(\gamma(t))=\lim_{s\rightarrow+\infty}g(\tilde{\gamma
}(s))-\liminf_{s\rightarrow+\infty}f(\tilde{\gamma}(s))=0.
\]
In either case we get $g(x_{0})\geq f(x_{0})$ and (since $x_{0}$ is arbitrary)
deduce $g\geq f$. Exchanging the roles of $g$ and $f$, we obtain the desired
result. \hfill$\Box$

\bigskip

\begin{remark}
[importance of integrability condition]\label{rem:tap} It is tempting to simplify Definition~\ref{def:salas} 
and replace the integrability condition in~(\ref{eq:jpp}) on the path $\gamma:[0,+\infty)\rightarrow\mathbb{R}^{d}$ by its consequence 
\begin{equation}
\underset{s\rightarrow+\infty}{\lim}T[g](\gamma(s))=0 \label{eq:tap0}
\end{equation}
(see (1) after Definition~\ref{def:salas} above), possibly together with the requirement that the curve is unbounded with no $\omega$-limits (to avoid reduncandy with critical points). The following example illustrates that this would lead to an undesirably large class of asymptotically critical curves: \medskip\newline
Indeed, set $\mathcal{U}
:=\mathbb{R}\times(0,+\infty)$ and consider the convex function
\begin{equation}
\left\{
\begin{array}
[c]{l}
g:\mathcal{U}\rightarrow\mathbb{R}\medskip\\
g(x,y)=\frac{x^{2}}{y}.
\end{array}
\right.  \label{eq:tap}
\end{equation}
Then, $\mathrm{Im}(g)= \lbrack0,+\infty)$ and $c=0$ is the unique critical value
of $g$ (in particular, every point of the level set $[g=0]:=\{0\}\times(0,+\infty)$ is critical\,!).
Since $g$ is a convex and $\inf g =0$, the result of~\cite{BCD2018} applies and $g$ is determined by its slope $\|\nabla g\|$ and its minimum value $\inf g=0$. Let us now focus on assumption (iii) of Theorem~\ref{thm:Baillon} (control on asymptotic critical paths). For a convex function, the integrability condition~(\ref{eq:jpp}) yields $g(\gamma(s))\rightarrow\inf g=0$. Therefore, (iii) is not an additional requirement (it is already contained in assumption~(ii)) and Theorem~\ref{thm:Baillon} generalizes the convex determination result mentioned in the introduction.

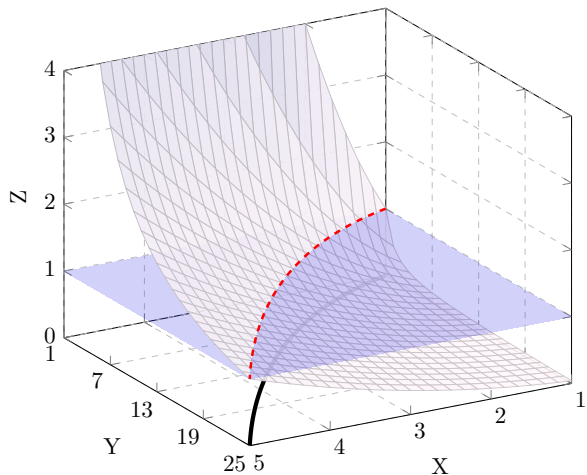
\begin{figure}[!h]
    \centering
    \begin{tikzpicture}[scale = 0.8]
    \begin{axis}[view={150}{25},
        grid = major,
        grid style={dashed, gray!50},
        scale only axis,
        xmin=1, xmax=5,
xtick={1,2,3,4,5},
xlabel={X},
ymin=1, ymax=25,
ytick={1,7,13,19,25},
ylabel={Y},
zmin=0, zmax=4,
zlabel={Z},
    ]

\addplot3[
        domain=1:5,
        samples = 60,
        samples y = 0,
        line width=2pt
    ] (
{x},
{x^2},
{0}
);

\addplot3[surf, mesh/cols=3, shader=interp, opacity = 0.2]  table [col sep=semicolon,x=x,y=y,z=z] {PlaneLow.txt}; 
    
    \addplot3 [surf, domain= 1:5,
    colormap/PuBu,
    y domain=1:25, opacity=0.9]
    {x^2/y};

\addplot3[surf, mesh/cols=3, shader=interp, opacity = 0.2]  table [col sep=semicolon,x=x,y=y,z=z] {PlaneUp.txt};

\addplot3[domain=1:5,
        red,
        very thick, dashed,
        samples=60,
        samples y=0,
    ] (
{x},
{x^2},
{1}
); 

    \end{axis}
    \end{tikzpicture} 
    \caption{{\footnotesize Function $g(x,y) = x^2/y$. In blue, the plane $z=1$.  In black, the curve $\gamma(t)$ for $c = 1$. In dashed red, the curve $(\gamma(t),g(\gamma(t)))=(\gamma(t),c)$ for $c=1$. Plane XY plotted in 1:6, starting at point (1,1). Plane XZ plotted in~1:1, starting at point (1,0).}}
    \label{fig:IntegrabilityExample}
\end{figure}

Omiting the integrability condition in Definition~\ref{def:salas} would have led to a completely different situation: for every $c>0$, the level set 
$$[g=c]:=\{(x,y)\in\mathcal{U}: x^2=cy\}=\{(\sqrt{c}t,t^2):t\neq 0 \}$$
contains an unbounded curve $\gamma(t)=(\sqrt{c}t,t^{2})$, $t>0$, without $\omega$-limits,
satisfying
$$g(\gamma(t))=c\qquad \text{and }\qquad\nabla g(\gamma(t))=\left( \frac{2\sqrt{c}} {t},\frac{c}{t^{2}}\right)  \underset{t\rightarrow\infty}{\longrightarrow}0,$$
that is, every $c>0$ would have been an asymptotic critical value of $g$ (see Figure~\ref{fig:IntegrabilityExample} for an illustration of the case $c=1$). Therefore, assumption~(iii) of Theorem~\ref{thm:Baillon} becomes very restrictive leading to an essentially useless statement.\hfill$\Diamond$
\end{remark}
\medskip

In what follows, we will move to metric spaces and establish determination
results for general classes of functions. The lack of derivatives (and norms)
is addressed by the metric slope, which we consider under an abstract unified
framework encompassing several other paradigms of descent operators. An
additional difficulty is to control the asymptotic behaviour of the functions:
asymptotically critical paths will be replaced by asymptotically critical
sequences $\{z_{n}\}_{n}$ with $T[g](z_{n})\rightarrow0$ as $n\rightarrow
\infty$ and the integrability condition~(\ref{eq:jpp}) by the summability
condition
\begin{equation}
{\displaystyle\sum\limits_{n\geq1}}\,T[g](z_{n})\,d(z_{n},z_{n+1}
)\,<\,+\infty\text{.} \label{eq:kruger}
\end{equation}
Indeed, since $\tilde{\gamma}$ is parameterized by arc-length, setting
$z_{n}:=\tilde{\gamma}(s_{n}),$ for all $n\geq1,$ we deduce
\[
s_{n+1}-s_{n}\geq d(\tilde{\gamma}(s_{n+1}),\tilde{\gamma}(s_{n}))
\]
and (\ref{eq:kruger}) follows directly from the discretization
of~(\ref{eq:jpp}). This together with the fact that the sequence has no
accumulation points eventually leads to Definition~\ref{def:AC-1}.


\subsection{General descent paradigms}

\label{ssec:2.3}Our goal is to derive a nonsmooth version of
Theorem~\ref{thm:Baillon} for functions defined in a complete metric space
$(X,d)$. As already mentioned, this requires a suitable extension of the
notion of asymptotically critical paths in order to impose a boundary
condition in the lines of condition~(iii) of Theorem~\ref{thm:Baillon}.
Moreover, we aim to obtain a statement that generalizes (and recovers) the
determination results of \cite{DS2022} (local slope, inf-compact case) and
\cite{TZ2022} (global slope, complete metric case) and apply to the main
paradigms of descent moduli studied in \cite{DMS2022} and quoted below:

\begin{enumerate}
\item \textit{The local (metric) slope} (\cite{GMT1980}, \cite{AGS2008}
\emph{e.g.}): For any $f\in\mathcal{LSC}(X)$ the local (metric) slope is given by:

\begin{equation}
s[f](x)=s_{f}(x):=\left\{
\begin{array}
[c]{cl}
\displaystyle\limsup_{y\rightarrow x}\,\,\Delta^{+}f(x,y), & \text{ if }
x\in\dom f\\
+\infty, & \text{ otherwise.}
\end{array}
\right.  \label{eq:metricSlope}
\end{equation}

\item \textit{The global slope}: Similarly, for any $f\in\mathcal{LSC}(X)$ the
global slope is given by

\begin{equation}
\mathcal{G}[f](x):=\left\{
\begin{array}
[c]{cl}
\underset{y\in X}{\sup}\,\,\Delta^{+}f(x,y), & \text{ if }x\in\dom f\\
+\infty, & \text{ otherwise.}
\end{array}
\right.  \label{eq:globalSlope}
\end{equation}

\item \textit{The average descent modulus}: Let $\mu$ be a probability measure
over $X$, and let $\mathcal{F}$ be the $\mu$-measurable extended-valued
functions on $X$. The average descent modulus is then given by
\begin{equation}
\mathcal{M}[f](x)=\left\{
\begin{array}
[c]{cl}
\displaystyle\int_{X}\Delta_{f}^{+}(x,y)\mu(dy)=\int_{X}\left\{
f(x)-f(y)\right\}  ^{+}\left(  \frac{1}{d(x,y)}\mu(dy)\right)  , & \text{if
}x\in\dom f,\\
& \\
+\infty, & \text{otherwise}.
\end{array}
\right.  \label{eq:nonlocalOp}
\end{equation}
This is an oriented nonlocal operator determined by the family of measures
$\{\mu_{x}\}_{x}$ with
\[
\mu_{x}(dy)=\frac{1}{d(x,y)}\mu(dy),\,\,\,\text{if }y\neq x\qquad\text{(under
the convention }\frac{0}{0}=0\text{).}
\]
(see \cite[Definition~4.14]{DMS2022}).

\item \textit{The diffusion descent modulus}: Let $\mathcal{F}$ be the $\mu
$-measurable extended-valued functions on $X$, and suppose now that $\mu(A)>0$
for every open set $A$ of $X$. Then, we define the diffusion descent modulus
$\mathcal{D}$ over $\mathcal{F}$ given by
\begin{equation}
\mathcal{D}[f](x)=\left\{
\begin{array}
[c]{cl}
\displaystyle\limsup_{\varepsilon\rightarrow0^{+}}\frac{1}{\mu(B(x,\varepsilon
))}\int\limits_{B(x,\varepsilon)}\Delta_{f}^{+}(x,y)\mu(dy), & \text{ if }
x\in\dom f,\\
& \\
+\infty, & \text{otherwise.}
\end{array}
\right.  \label{eq:DispersionOp}
\end{equation}
This is the oriented 1-dispersion operator for measure $\mu$ (see
\cite[Definition 4.2]{DMS2022}).
\end{enumerate}

All four notions described above fit the definition of descent modulus in the
extended sense of Definition~\ref{def:T}. The proofs are mild adaptations of
\cite{DMS2022}. In all cases, Lipschitz continuous functions are contained in
the domain of each of the aformentioned descent moduli.


\subsection{Metrically compatible descent moduli}

\label{ssec:2.4} The definition of a descent modulus (\emph{cf}.
Definition~\ref{def:T}) does not require prior assumptions on the space~$X$.
In particular, $X$ does not need to be metric (neither topological) space,
although the aforementioned determination result in \cite[Theorem~3.5]
{DMS2022} will eventually require topology, to formulate the assumptions of
continuity and compactness. This being said, this work is inscribed in the
framework of a complete metric space $(X,d).$ In order to obtain determination
results in this setting and ensure an efficient use of completeness property,
we need to impose some (metric) compatibility condition to the considered
descent modulus reckoning with completeness of $X$. The condition should
encompass the four paradigms of Subsection~\ref{ssec:2.3}. The first natural
attempt for such condition leads to the following definition.

\begin{definition}
[strong metric compatibility]\label{def:Tmetric_stark}We say that a descent
modulus $T:\mathcal{F}\rightarrow\lbrack0,+\infty]^{X}$ is strongly metrically
compatible, if for some strictly increasing continuous function $\theta
:\mathbb{R}_{+}\rightarrow\mathbb{R}_{+}$ with $\theta(0)=0$ and
$\underset{t\rightarrow\infty}{\lim}\theta(t)=+\infty$ and for every
$f,g\in\mathrm{dom}(T)$, $x\in\dom g$ and $\delta>0$ it holds:
\begin{equation}
T[f](x)<\delta<T[g](x)\quad\implies\quad\exists z\in\dom g:\,\;\dfrac
{\{f(x)-f(z)\}^+}{d(x,z)}\,<\,\theta(\delta)\,<\,\dfrac{\{g(x)-g(z)\}^+}{d(x,z)}
\tag{$\widehat{\mathcal{D}_2}$}
\end{equation}

\end{definition}

The idea behind the above definition is to guarantee that whenever $T[f](x)<\delta<T[g](x)$, there exists a point $z$ ensuring on the one hand more descent for $g$ than for $f$ and on
the other hand enough descent for $g$ relative to the distance $d(x,z)$ (up to
a factor $\theta(\delta)^{-1}$). \smallskip\newline
It is straightforward to see that $(\widehat{\mathcal{D}_{2}})$
yields $(\mathcal{D}_{2})$. 
Another important remark is that the function
$\theta$ is invertible, and its inverse $\theta^{-1}:\mathbb{R}_{+}
\rightarrow\mathbb{R_{+}}$ verify the same properties as $\theta$, that is,
$\theta$ is a strictly increasing continuous function with $\theta^{-1}(0)=0$
and $\underset{s\rightarrow\infty}{\lim}\theta^{-1}(s)=+\infty$.\smallskip
\newline The above definition encompasses the \textit{sup}--type paradigms of
metric descent modulus.

\begin{example}
[steepest descent operators]\label{ejem32}(i). (local slope) Let us recall
from~\eqref{eq:metricSlope} the definition of the local slope. It has been
shown in \cite[Proposition 3.7]{DMS2022} that the ``local slope'' operator
\[
s:(\mathbb{R}\cup\{+\infty\})^{X}\rightarrow\lbrack0,+\infty]^{X},
\]
defined by $s[f]:=s_{f}$ is a descent modulus on the linear subspace
$\mathcal{C}(X)$, that is, it satisfies properties $(\mathcal{D}_{1}
)$--$(\mathcal{D}_{3})$ of Definition~\ref{def:T}. Since $(\mathcal{D}_{0})$
is also verified, it follows easily that $s[\cdot]$ is a descent modulus on the
translation cone $\mathcal{F}:=\mathcal{LSC}(X)$ as well. Moreover, let
$f,g\in\mathcal{LSC}(X)$ and $x\in\dom g$ be such that $s_{f}(x)<\delta
<s_{g}(x).$ Then,~\eqref{eq:metricSlope} yields that for some $\sigma>0$
sufficiently small, we have:
\[
\sup_{y\in B(x,\sigma)}\Delta_{f}^{+}(x,y)<\delta<\sup_{y\in B(x,\sigma)}\Delta_{g}^{+}(x,y).
\]
In particular, we can choose $z\in B(x,\sigma)$ such that $\Delta_{g}^{+}(x,z)>\delta$, and obtain
\[
\dfrac{\{f(x)-f(z)\}^{+}}{d(x,z)}<\delta<\dfrac{\{g(x)-g(z)\}^{+}}{d(x,z)}.
\]
Thus $(\widehat{\mathcal{D}_{2}})$ holds for $\theta(\delta)=\delta
$.\smallskip\newline

\noindent(ii). (global slope) The global slope of a function $f\in\mathcal{LSC}(X)$ at a point $x\in X$ is defined in~\eqref{eq:globalSlope}
as follows:
\[
\mathcal{G}_{f}(x):=\,\sup_{y\in X}\Delta_{f}^{+}(x,y)=\sup_{y\in\dom f}
\Delta_{f}^{+}(x,y)
\]
It is straightforward to see that the \textquotedblleft global
slope\textquotedblright\ operator $\mathcal{G}:\mathcal{LSC}(X)\rightarrow
\lbrack0,+\infty]^{X}$ is strongly metrically compatible, satisfying
$(\widehat{\mathcal{D}_{2}})$ with $\theta(\delta)=\delta$.\hfill$\Diamond$
\end{example}

\bigskip

Notice that if a descent modulus $T$ is strongly metrically
compatible, then $\hat{T}:=\phi\circ T$ remains strongly metrically
compatible, for every strictly increasing continuous function with $\phi(0)=0$
and $\underset{t\rightarrow+\infty}{\lim}\phi(t)=+\infty$. 
However, unfortunately, Definition~\ref{def:Tmetric_stark} is quite restrictive and
fails to cope with some important \textit{average}--type paradigms, as reveals
the following example:


\begin{example}[Average descent fails strong metric compatibility]\label{ex:AverageNotStronglyMetricCompatible}
Set $X =[0,+\infty)$ endowed with the distance $d(t,s) = |t-s|$ and the usual Lebesgue measure. We show that for any function $\theta:[0,+\infty)\to [0,+\infty)$ given as in Definition~\ref{def:Tmetric_stark} and any $\delta>0$, there exist functions $f$ and $g$ for which ($\widehat{\mathcal{D}_2}$) fails for the average descent modulus $\mathcal{M}$ (see Subsection~\ref{ssec:2.3}). Indeed, fix $\theta$ and $\delta>0$ and set 
\[
t_0 = \frac{3\delta}{5\theta(\delta)}.
\]
 
Let us define functions $\psi,\phi:X\to [0,+\infty)$ as follows:
\[
\psi(t) = \left\{\begin{array}{cl}
    \theta(\delta)&\text{ if }t\in[0,t_0)\smallskip \\
     \frac{3}{2}\theta(\delta) - \frac{5\theta(\delta)^2}{6\delta}t& \text{ if }t\in [t_0,3t_0], \smallskip\\
     0& \text{ if }t>3t_0. 
\end{array}\right.\quad\text{and}\quad \phi(t) = \left\{\begin{array}{cl}
        \theta(\delta)&\text{ if }t\in[0,t_0),\smallskip \\
     2\theta(\delta) - \frac{5\theta(\delta)^2}{3\delta}t& \text{ if }t\in [t_0,2t_0],\smallskip \\
     0& \text{ if }t>2t_0. 
\end{array}\right.
\]
It is not hard to realize (see Figure~\ref{fig:dg-and-df}) that
\[
\int_X \psi(t) dt = \frac{6}{5}\delta \, > \, \delta \, > \, \frac{9}{10}\delta = \int_X \phi(t)dt.
\]
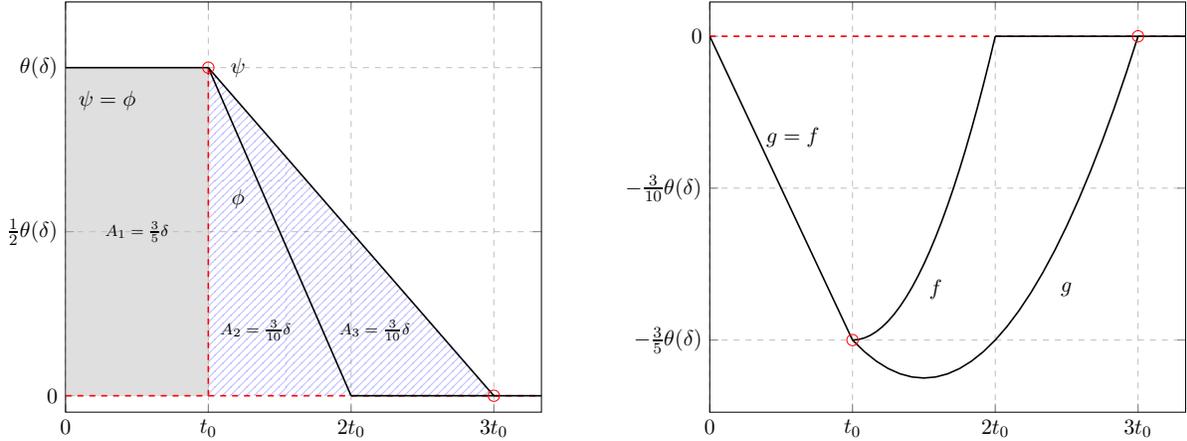
\begin{figure}[!h]
\begin{minipage}{0.5\textwidth}
    \begin{tikzpicture}[scale=0.75]
    \begin{axis}[
        grid = major,
        grid style={dashed, gray!50},
        scale only axis,
        xmin=0,xmax=4,
        ymin = -0.025, ymax= 0.6,
        xtick={0,1.2,2.4,3.6},
        xticklabels={0,$t_0$,$2t_0$,$3t_0$},
        ytick={0,0.25,0.5},
        yticklabels={0,$\tfrac{1}{2}\theta(\delta)$,$\theta(\delta)$},
    ]
        \filldraw[pattern=north east lines, pattern color=blue,draw=none,opacity = 0.5] (axis cs: {2.4},{0})--(axis cs: {3.6},{0}) -- (axis cs: {1.2},{0.5})--cycle;
        \filldraw[pattern=north east lines, pattern color=blue,draw=none,opacity = 0.5] (axis cs: {1.2},{0})--(axis cs: {1.2},{0.5}) -- (axis cs: {2.4},{0})--cycle;
        \fill[color = gray!50!white, opacity=0.5] (axis cs:{0},{0}) -- (axis cs:{1.2},{0}) -- (axis cs:{1.2},{0.5}) -- (axis cs:{0},{0.5})--cycle;
        
        \addplot[red,thick,dashed,domain=0:4]{0};
        \addplot[black,thick, domain = 1.2:2.4]{ 1 - (5/12)*x };
        \addplot[black,thick, domain = 1.2:3.6]{ 0.75 - (5*0.25/6)*x };
        \addplot[black,thick, domain = 2.4:4]{ 0 };
        \addplot[black,thick, domain = 0:1.2]{ 0.5 };
        \draw[red,thick,dashed] (axis cs:{1.2},{0})--(axis cs:{1.2},{0.5});
        \draw[red] (axis cs: {1.2},{0.5}) circle (0.1cm);
        \draw[red] (axis cs: {3.6},{0}) circle (0.1cm);

        \node at (axis cs: {0.35},{0.45}){$\psi=\phi$};
        \node at (axis cs: {1.45},{0.5}){$\psi$};
        \node at (axis cs: {1.45},{0.3}){$\phi$};
        \node at (axis cs: {0.6},{0.25}){\scriptsize $A_1= \frac{3}{5}\delta$};
        \node at (axis cs: {1.6},{0.1}){\scriptsize $A_2 = \frac{3}{10}\delta$};
        \node at (axis cs: {2.6},{0.1}){\scriptsize $A_3 = \frac{3}{10}\delta$};
    \end{axis}
    \end{tikzpicture}
    \end{minipage}
\begin{minipage}{0.45\textwidth}
   \centering
    \begin{tikzpicture}[scale=0.75]
    \begin{axis}[
        grid = major,
        grid style={dashed, gray!50},
        scale only axis,
        xmin=0,xmax=4,
        xtick={0,1.2,2.4,3.6},
        xticklabels={0,$t_0$,$2t_0$,$3t_0$},
        ytick={0,-0.3,-0.6},
        yticklabels={0,$-\frac{3}{10}\theta(\delta)$,$-\frac{3}{5}\theta(\delta)$},
    ]
     \addplot[red,thick,dashed,domain=0:4]{0}; 
     \addplot[black,thick, domain=1.2:2.4]{-(1 - (5/12)*x)*x};
     \addplot[black,thick, domain=2.4:4]{0};
     \addplot[black,thick, domain=0:1.2]{-0.5*x};
     \addplot[black,thick, domain=1.2:3.6]{-( 0.75 - (5*0.25/6)*x)*x};
     \draw[red] (axis cs: {1.2},{-0.6}) circle (0.1cm);
     \draw[red] (axis cs: {3.6},{0}) circle (0.1cm);

        \node at (axis cs: {3.0},{-0.5}){$g$};
        \node at (axis cs: {0.7},{-0.2}){$g=f$};
        \node at (axis cs: {1.9},{-0.5}){$f$};
    \end{axis}
    \end{tikzpicture}
    \end{minipage}
    \caption{ {\footnotesize \textit{Left:} Regions below the functions $\psi$ and $\phi$. \textit{Right:} The functions $g(t)=-t\psi(t)$ and $f(t) = -t\phi(t)$ coincide on $[0, t_0)$ while for $t\geq t_0$ we have $f\geq g$. }\label{fig:dg-and-df}}
\end{figure}

Now, we can define the functions $g,f:X\to \R$ given by 
\[
g(t) = -t\,\psi(t)\qquad \text{ and }\qquad f(t)= -t\,\phi(t).  
\]
Then $\psi$ and $\phi$ coincide with the functions $\Delta^+g(0,\cdot)$ and $\Delta^+f(0,\cdot)$. Indeed, since $g$ and $f$ attain their global maximum at $t=0$, we can directly write
\[
\Delta^+g(0,t) = \frac{g(0) - g(t)}{t} = \psi(t)\quad\text{and}\quad\Delta^+f(0,t) = \frac{f(0) - f(t)}{t} = \phi(t).
\]
Thus, ($\widehat{\mathcal{D}_2}$) fails at $0$, since $\mathcal{M}[f](0) = \frac{9}{10}\delta < \delta<\frac{6}{5}\delta <\mathcal{M}[g](0)$ but there is no $t\in X$ such that $\phi(t)< \theta(\delta) < \psi(t)$. Even worse, there is no $t\in X$ such that $\phi(t)<\psi(t)$ (more descent for $g$ than for $f$) and $\theta(\delta)\leq\psi(t)$  (enough descent for $g$) simultaneously.\hfill$\Diamond$
\end{example}

\bigskip

The previous example reveals that average--type descents fail Definition~\ref{def:Tmetric_stark} (strong metric compatibility). To overcome this difficulty, we
need to relax this definition as follows:

\begin{definition}
[metric compatibility]\label{def:Tmetric}A descent modulus $T:\mathcal{F}
\rightarrow\lbrack0,+\infty]^{X}$ is said to be metrically compatible, if for
every $\rho>0$, there exists a strictly increasing continuous function
$\theta_{\rho}:\mathbb{R}_{+}\rightarrow\mathbb{R}_{+}$ with $\theta_{\rho
}(0)=0$ and $\underset{t\rightarrow+\infty}{\lim}\theta_{\rho}(t)=+\infty$,
such that for every $f,g\in\mathrm{dom}(T)$, $x\in\dom g$ and $\delta>0$ it
holds:
\begin{equation}
T[f](x)<\delta<T[g](x)\,\,\implies\quad\exists z\in\dom g:\,\;\left\{
\begin{array}
[c]{c}
\{f(x)-f(z)\}^{+}<(1+\rho)\{g(x)-g(z)\}^{+}\smallskip \\
\text{and}\\
\theta_{\rho}(\delta)\,d(x,z)\,<\,g(x)-g(z)
\end{array}
\right.  \tag{$C$}
\end{equation}

\end{definition}

The above definition is a trade-off between the needs of the proof of
Theorem~\ref{thm.strict_compa} and a common scheme that incorporates all main
paradigms of Subsection~\ref{ssec:2.3}. The difference with
Definition~\ref{def:Tmetric_stark} is that given $\rho>0,$ if $T[f](x)<\delta
<T[g](x)$, we can find a point $z$ (depending on $\rho$) ensuring more descent
for $g$ than for $f$ up to a factor $(1+\rho)$ and sufficient descent for $g$
relative to the distance $d(x,z)$, up to a factor $\theta_{\rho}(\delta)^{-1}$
(depending again on $\rho$). In this sense, for a given tolerance $\rho>0$,
condition ($C$) is a trade-off between these requirements.\smallskip
\newline\textit{Average descent} operators are operators for which the descent
of $f$ at a point $x$ is obtained by integrating the quotient $\Delta_{f}^{+}(x,y)$ with respect to some probability measure on $X.$ This category of
operators are now metrically compatible with respect to this relaxed
definition:\smallskip

\begin{example}
[average descent moduli]\label{ejem26}Let $\mu$ be a probability measure on
the metric space $(X,d)$ and consider the translation cone
\[
\mathcal{F}:=\{f:X\rightarrow\mathbb{R}\cup\{+\infty\}\ :\ f\text{ is proper
and }\mu\text{-measurable}\}.
\]

\noindent(i) Consider the operator $\mathcal{M}$ defined in
(\ref{eq:nonlocalOp}). Let us show that $\mathcal{M}$ is a metrically
compatible descent modulus. Indeed, fix $\rho>0$ and assume that for some
$f,g\in\mathrm{dom}(\mathcal{M})$, $x\in\dom g$ and $\delta>0$ we have
$\mathcal{M}[f](x)<\delta<\mathcal{M}[g](x).$ Then, by definition, we have
that
\[
\int_{X}(1+\rho)\,\Delta_{g}^{+}(x,y)\,\mu(dy)\,>\,\int_{X}(\Delta_{f}^{+}(x,y)+\rho\,\delta)\,\mu(dy)\,=\,\int_{X}\Delta_{f}^{+}(x,y)\,d\mu
\,\,+\,\rho\,\delta,
\]
which yields that there exists $z\in\dom g$ such that $(1+\rho)\Delta_{g}^{+}(x,z)>\Delta_{f}^{+}(x,z)+\rho\delta$. The above inequality yields
\[
(1+\rho)\{g(x)-g(z)\}^{+}>\{f(x)-f(z)\}^{+}\quad\text{ and }\quad\frac
{g(x)-g(z)}{d(x,z)}>\underbrace{\left(  \frac{\rho}{1+\rho}\right)  \,\delta
}_{\theta_{\rho}(\delta)}.
\]
Thus, $(C)$ is verified for $\theta_{\rho}(\delta):=\frac{\rho}{1+\rho}\delta
$.\smallskip\newline

\noindent(ii) Suppose now that $\mu(A)>0$ for every open set $A$ of $X$ and
consider the oriented 1-dispersion operator for the measure $\mu,$ given by
(\ref{eq:DispersionOp}). Then, $\mathcal{D}$ is a metrically compatible descent
modulus. Indeed, let $f,g\in\mathrm{dom}(\mathcal{D})$, $x\in\dom g$ and
$\delta\in\mathbb{R}$ such that $\mathcal{D}[f](x)<\delta<\mathcal{D}[g](x)$.
Then, for $\varepsilon>0$ sufficiently small we have
\[
\frac{1}{\mu(B(x,\varepsilon))}\int\limits_{B(x,\varepsilon)}\Delta_{f}
^{+}(x,y)\mu(dy)\,<\,\delta<\,\frac{1}{\mu(B(x,\varepsilon))}\int%
\limits_{B(x,\varepsilon)}\Delta_{g}^{+}(x,y)\mu(dy).
\]
Then, the conclusion follows by noting that $\nu(dy)=\frac{1}{\mu
(B(x,\varepsilon))}\mu(dy)$ is a probability measure over the metric space
$B(x,\varepsilon)$, and therefore we can proceed as in the previous example
(i).\hfill$\Diamond$
\end{example}

\noindent Based on these examples-schemes, we can significantly enlarge the
class of metrically comptatible descent moduli as follows: for any strictly
increasing, continuous function $\phi:\mathbb{R}_{+}\rightarrow\mathbb{R}_{+}$
with $\phi(0)=0$ and $\underset{t\rightarrow+\infty}{\lim}\phi(t)=+\infty$, we
can replace $\Delta_{f}^{+}(x,y)$ (in definitions~\eqref{eq:metricSlope},
\eqref{eq:globalSlope}, \eqref{eq:nonlocalOp} and \eqref{eq:DispersionOp}) by
\[
\widetilde{\Delta_{f}^{+}}(x,y):=\phi^{-1}(\Delta_{f}^{+}(x,y))
\]
and obtain new classes.\smallskip\newline

Let us now give an example of a descent modulus which is not metric compatible.

\begin{example}
Take $X=\mathbb{N}$ endowed with the distance function $d$ given by
$d(m,n)=|m-n|$. Consider the operator $T$ given by
\[
T[f](k)=\sup_{m\in\mathbb{N}}\,\left\{  f(k)-f(m)\right\}  ^{+},\qquad
\text{for all } k\in\mathbb{N}.
\]
Clearly, $T$ is a descent modulus (it coincides with the global slope for the
distance $d_{0}$ given by $d_{0}(k,m)=1$, if $k\neq m$, and $d_{0} (k,m)=0$,
if $k=m$). Consider the function $f_{n}$, $n\in\mathbb{N}$, given by
\[
f_{n}(m)=
\begin{cases}
1, & \text{ if }m\neq n,\\
0, & \text{ if }m=n.
\end{cases}
\]
Observe that for each $\delta\in(0,1)$ and each $n\geq2$, we have that
$T[f_{n}](1)=1>\delta>0=T[f_{1}](1)$. However, regardless the tolerance
$\rho>0$, the only choice for $m\in\lbrack f_{n}<f_{n}(1)]$ is $m=n$ and so
\[
(1+\rho)\,\Delta^{+}_{f_{n}}(1,n)\,=\,(1+\rho)\,\frac{f_{n}(1)-f_{n}
(n)}{d(1,n)} =\frac{1+\rho}{n-1} \,>\,0=\frac{\left\{  f_{1}(1)-f_{1}
(n)\right\}  ^{+}}{d(1,n)}=\Delta^{+}_{f_{1} }(1,n).
\]
If $T$ were metrically compatible, there would exist a continuous, strictly
increasing function $\theta_{\rho}:[0,+\infty)\rightarrow\lbrack0,+\infty)$,
with $\theta_{\rho}(0)=0$ and $\theta_{\rho}(\delta)<\frac{f_{n}(1)-f_{n}
(n)}{d(n,1)}=\frac{1}{n-1}$ for all $\delta\in(0,1)$ and all $n\geq2,$ which
is a contradiction. Therefore $T$ is not metrically compatible for the metric
$d$. Notice however, that it is so for the discrete metric $d_{0}.$
\hfill$\Diamond$
\end{example}

\begin{remark}
We recall from \cite[Section 3.4]{DMS2022} that there exist slope--like
operators, as the weak slope (\cite{CDM1993}, \cite{DM1994Critical} e.g.) or
the limiting slope (\cite{DLI2015}) that are not descent moduli, since they
fail property $(\mathcal{D}_{2})$ (monotonicity).
\end{remark}

\noindent We finish the section with two stability properties that we will
need in the sequel. The first one stems from \cite[Proposition~3.2(b)]{DMS2022}, 
whose proof is based on the monotonicity property $(\mathcal{D}_{2})$ of the descent moduli.

\begin{proposition}
[translation-invariance]\label{prop.mono_proper}Let $T:\mathcal{F}
\rightarrow\lbrack0,+\infty]^{X}$ be a descent modulus. Then for any
$f\in\mathcal{F}$ and $c\in\mathbb{R}$, we have
\[
T[f]=T[f+c].
\]

\end{proposition}

The second property provides another way of constructing new metrically
compatible descent moduli, similar to the remark after Example~\ref{ejem26}.

\begin{proposition}
[composition with increasing functional]Let $\phi:\mathbb{R}_{+}
\rightarrow\mathbb{R}_{+}$ be a strictly increasing continuous function with
$\phi(0)=0$ and $\underset{t\to+\infty}{\lim}\phi(t)=+\infty$. Let
\[
T:\mathcal{F}\rightarrow\lbrack0,+\infty\rbrack^{X}
\]
be a descent modulus. One has that
\begin{equation}
T\text{ is metrically compatible }\iff\phi\circ T\text{ is metrically
compatible}, \label{eq:Composition}
\end{equation}
where $\phi\circ T$ is the descent modulus given by $(\phi\circ T)[f](x)=\phi
(T[f](x))$, with $\phi(+\infty)=+\infty$.
\end{proposition}

\noindent\textbf{Proof.} The fact that $\phi\circ T$ is a descent modulus was
established in \cite[Proposition~3.9]{DMS2022} for real-valued functions. The
proof can be easily adapted to the present setting of extended real-valued
functions. Concerning metric compatibility, since $\phi^{-1}$ is also strictly
increasing and continuous, with $\phi^{-1}(0)=0$ and $\underset{s\to
+\infty}{\lim}\phi^{-1}(s)=+\infty$, it is sufficient to establish only one
implication. To this end, let us assume that $T$ is metrically compatible with
$\{\theta_{\rho}\}_{\rho>0}$ given as in Definition~\ref{def:Tmetric}. Let
$f,g\in\mathrm{dom}(\phi T)$, $x\in\dom g$ and $\delta>0$ such that $(\phi
T)[f](x)<\delta<(\phi T)[g](x)$. Fix $\rho>0$. It is straightforward to see
that $f,g\in\mathrm{dom}(T)$ and $T[f](x)<\phi^{-1}(\delta)<T[g](x)$,
therefore for some $z\in\dom g$ we have
\[
\{f(x)-f(z)\}^{+} < (1+\rho)\{g(x) - g(z)\}^{+}\quad\text{ and }\quad
\theta_{\rho}(\phi^{-1}(\delta))<\frac{g(x)-g(z)}{d(x,z)}.
\]
Since $\widetilde{\theta_{\rho}}:=\theta_{\rho}\circ\phi^{-1}$ is continuous
and strictly increasing with $\widetilde{\theta_{\rho}}(0)=0$ and
$\displaystyle\lim_{t\rightarrow+\infty}\widetilde{\theta_{\rho}}(t)=+\infty$,
the conclusion follows. \hfill$\Box$

\section{Main result}
\label{sec:3}

This section contains our main result: we establish a comparison principle, in
the lines of Lemma~\ref{lemma_sal1}, for metrically compatible descent moduli
(\emph{cf.} Definition~\ref{def:Tmetric}) in a complete, but not necessarily
compact, metric space. This result will eventually lead to our determination
result (Theorem~\ref{thm:DetExtended-metric}). Absence of compactness imposes
an assumption on the asymptotic behaviour for which, as already discussed in
Subsection~\ref{ssec:2.2}, the choice of the notion of asymptotic criticality
is paramount. This latter not only consists of saying that the descent moduli
vanish at infinity, but also requires two additional restrictions: absence of
accumulation points and a summability condition (Definition~\ref{def:AC-1}).
The price to pay is that the scheme
of proof becomes more involved, but as a reward, we are able to obtain a
statement that generalizes all previous results (\emph{c.f.}
Theorem~\ref{thm:DetExtended-metric}). Both the local slope determination of
\cite{DS2022} for inf-compact functions and the global slope determination
\cite{TZ2022} are now recovered by our final statement.

\subsection{Comparison lemmas in complete metric spaces}

\label{sec:Strict}

In this section, we want to establish a comparison principle, in the lines of
Lemma~\ref{lemma_sal1}, for metrically compatible descent moduli (\emph{cf.}
Definition~\ref{def:Tmetric}) in a complete (but not necessarily compact)
metric space. Absence of compactness assumption leads inevitably to impose
control on the asymptotic behavior. To do so, we need a precise notion of
asymptotic $T$-criticality, which is given in the following definition:

\begin{definition}
[asymptotic critical sequences]\label{def:AC-1}A sequence $\{z_{n}\}_{n\geq
1}\subset X\setminus\mathcal{Z}_{T}(g)$ is called $T$-asymptotically critical
for a function $g\in\mathcal{F}$ (in short, $T[g]$-critical) if it has no
converging subsequence and
\begin{equation}
\sum_{n=0}^{+\infty}
T[g](z_{n})\,d(z_{n},z_{n+1})<+\infty. \label{eq:AC}
\end{equation}
We denote by $\mathcal{AZ}_{T}(g)$ the set of asymptotically critical
sequences for $g$.
\end{definition}

\begin{remark}
[justification of terminology]\label{rem-critical}Any sequence $\{z_{n}\}_{n\geq1}$ satisfying ~$\underset{n\rightarrow\infty}{\lim\inf}\,T[g](z_{n})>0$ and~(\ref{eq:AC}) is necessarily Cauchy: indeed, assume that
for some $\delta>0$ and $N\geq1$ we have $T[g](z_{n})\geq\delta,$ for all
$n\geq N.$ Then for every $\varepsilon>0$, we take $n_{0}\geq N$ such that
\[
\sum_{i\geq N}
T[g](z_{i})\,d(z_{i},z_{i+1})<\varepsilon\,\delta,
\]
and for every $m>n\geq n_{0}$ we deduce
\[
d(z_{n},z_{m})\leq
{\displaystyle\sum\limits_{i=n}^{m-1}}
\,d(z_{i},z_{i+1})\,\leq\,\frac{1}{\delta}\,
{\displaystyle\sum\limits_{i=n}^{m-1}}
\,T[g](z_{i})\,d(z_{i},z_{i+1})<\varepsilon.
\]
Therefore, absence of convergent subsequences in Definition~\ref{def:AC-1}
ensures that in a complete metric space, every $T[g]$-critical sequence
$\{z_{n}\}_{n}$ should satisfy
\begin{equation}
\underset{n\rightarrow\infty}{\lim\inf}\,T[g](z_{n})=0. \label{eq:AC1}
\end{equation}
Another important feature of this notion is that it becomes vacuous under
compactness.
\end{remark}

The requirement for the convergence of the series in~\eqref{eq:AC} restricts
significantly the class of critical sequences. This was motivated by the
discussion in Subsection~\ref{ssec:2.2} and can be seen as a discrete version
of the integrability condition in~(\ref{eq:jpp}), which in turn was inspired
by the steepest descent gradient system. A similar condition was also employed in
the fundamental lemma~\cite[Lemma~4.2]{TZ2022} used to derive the
determination result in \cite{TZ2022}. Our next result
(Theorem~\ref{thm.strict_compa}) extends the metric determination result of \cite{TZ2022},
since it holds for any metrically compatible descent modulus (and not only for
the global slope). Indeed, in case of the global slope we have (see
forthcoming Lemma~\ref{lem:jortega}):

\begin{itemize}
\item If $T[g]=\mathcal{G}[g]$ (global slope) and $\{z_{n}\}_{n}$ is a
$T[g]$--critical sequence, then $g(z_{n})\underset{n\rightarrow\infty
}{\longrightarrow}\inf g$.
\end{itemize}

We are now ready to establish the main comparison lemma.

\begin{theorem}
[comparison lemma]\label{thm.strict_compa} Let $(X,d)$ a complete metric space
and $T:\mathcal{F}\rightarrow\lbrack0,+\infty]^{X}$ a metrically compatible
descent modulus. Let $f,g\in\mathrm{dom}(T)$ be two bounded from below
functions with $g\in\mathcal{LSC}(X)$ and $f\in\overline{\mathcal{C}}(X)$. Let
us assume that:

\begin{itemize}
\item[$\mathrm{(i).}$] (descent domination) $T[f](x)<T[g](x),$ for every $x\in
X\setminus\mathcal{Z}_{T}(g)$ ;

\item[$\mathrm{(ii).}$] (control of criticality) $f(z)\leq g(z),$ for every
$z\in\mathcal{Z}_{T}(g)$ ;

\item[$\mathrm{(iii).}$] For every $\{z_{n}\}_{n}\in\mathcal{AZ}_{T}(g)$ it
holds:
\begin{equation}
\underset{n\rightarrow+\infty}{\liminf}\mathrm{\,}f(z_{n})\,\leq
\,\underset{n\rightarrow+\infty}{\liminf}\mathrm{\,}g(z_{n}).
\label{eq:morduk}
\end{equation}

\end{itemize}

Then, it holds
\[
f(x)\leq g(x),\qquad\text{for every }x\in\dom g.
\]

\end{theorem}

\noindent\textbf{Proof.} We deduce from (i) and $(\mathcal{D}_{0})$ of
Definition~\ref{def:T} that $\dom g\subset\dom f$. \smallskip\newline Let us
fix $\rho>0$. Replacing $f$ and $g$ by $f-\inf g+1$ and $g-\inf g+1$ if
needed, we may assume that $g>0$ and consequently $(1+\rho)g>g>0$. In what
follows, we prove $f<(1+\rho)g$. Once this is done, since $\rho$ is
arbitrarily small, we deduce $f\leq g.$\smallskip\newline To this end, let
$\{\theta_{\rho^{\prime}}\}_{\rho^{\prime}>0}$ be a family of continuous
strictly increasing functions associated to the descent modulus $T$
(\textit{c.f.} Definition~\ref{def:Tmetric}). By
Proposition~\ref{eq:Composition} the operator $\widehat{T}:=\theta_{\rho}\circ
T$ is a metrically compatible descent modulus, which preserves $T$-critical
points (i.e. $\mathcal{Z}_{\widehat{T}}(g)=\mathcal{Z}_{T}(g)$) and
$T[g]$-critical sequences (i.e. a sequence is $T[g]$-critical if and only if
it is $\widehat{T}[g]$-critical). Furthermore, assumption (i) continues to
hold for $\widehat{T}$ and $(C)$ is now satisfied for $\theta_{\rho}
(\delta)=\delta$, $\delta>0$. Therefore, by replacing $T$ by $\widehat{T}$ if
necessary, we may assume that for the value $\rho>0$ (that we fixed in the
beginning) condition $(C)$ holds for the identity function $\theta_{\rho
}(\delta)=\delta$.\smallskip\newline With all the above in mind, let us
define, for every $x\in\dom g,$ the quantity
\begin{equation}
\delta(x):=\frac{1}{2}T[f](x)+\frac{1}{2}T[g](x). \label{eq:delta}
\end{equation}
Notice that since $f,g\in\dom T,$ we have $\delta(x)<+\infty.$ Moreover, if
$x\notin\mathcal{Z}_{T}(g)$ (that is, $x$ is not $T$-critical for $g$), then
(i) yields $\delta(x)>0$ and
\begin{equation}
T[f](x)<\delta(x)<T[g](x)<2\delta(x). \label{eq:delta1}
\end{equation}
We now assume, towards a contradiction, that there exists $x_{0}\in\dom g$
such that
\[
f(x_{0})\geq(1+\rho)\,g(x_{0}).
\]
In what follows, using successively assumption (i) (descent domination) and
condition $(C)$ of Definition~\ref{def:Tmetric}, we are going to construct a
sequence of points $\{x_{n}\}_{n\geq1}$ failing the conclusion of our theorem,
where additionally $g$ is strictly decreasing.

\medskip

\emph{Basic iteration scheme }(classical induction). Thanks to assumption (ii),
we have $x_{0}\notin\mathcal{Z}_{T}(g)$ and we deduce from (i) that
(\ref{eq:delta1}) holds for $x=x_{0}$. Using $(C)$ for $\theta(\delta)=\delta$
and $\delta=\delta(x_{0})$ we obtain $x_{1}\in\lbrack g<g(x_{0})]$ such that
\begin{equation}
0<f(x_{0})-(1+\rho)g(x_{0})\leq f(x_{1})-(1+\rho)g(x_{1})\quad\text{and\quad
}\delta(x_{0})\,d(x_{0},x_{1})\leq g(x_{0})-g(x_{1}).\label{eq:it1}
\end{equation}
Therefore we obtain:
\[
g(x_{0})>g(x_{1}),\qquad f(x_{1})>(1+\rho)g(x_{1})\qquad\text{and}\qquad
x_{1}\notin\mathcal{Z}_{T}(g).
\]
It is quite clear that the above procedure can be repeated as many times as we
wish, producing a sequence $\{x_{n}\}_{n\geq1}$ in $X\setminus\mathcal{Z}
_{T}(g)$ such that
\begin{equation}
\text{the sequence }\quad\left\{  \left(  f-(1+\rho)g\right)  (x_{n})\right\}
_{n\geq0}\text{ }\quad\text{is nondecreasing and positive}\label{eq:boris0}
\end{equation}
and
\begin{equation}
0\,<\,\delta(x_{n})\,d(x_{n},x_{n+1})\,<\,g(x_{n})-g(x_{n+1}),\text{ }
\quad\text{for every }n\geq0.\label{eq:boris1}
\end{equation}
The telescopic series obtained by summing up the above inequalities for all
$n\geq0$, together with (\ref{eq:delta1}) and the fact that the sequence
$\left\{  g(x_{n})\right\}  _{n\geq0}$ is strictly decreasing and bounded from
below, yield that
\begin{equation}
\frac{1}{2}\,{\displaystyle\sum\limits_{n=0}^{\infty}}
T[g](x_{n})\,d(x_{n},x_{n+1})\,<\,
{\displaystyle\sum\limits_{n=0}^{\infty}}
\delta(x_{n})\,d(x_{n},x_{n+1})\,\leq\,g(x_{0})-\inf g<+\infty
.\label{eq:fat-boy}
\end{equation}
Assumption (iii) together with (\ref{eq:boris0}) ensure that the sequence
$\{x_{n}\}_{n\geq0}$ cannot be $T[g]$--critical, therefore we deduce from the
above inequality and Definition~\ref{def:AC-1} that $\{x_{n}\}_{n\geq0}$ has
accumulation points as $n\rightarrow\infty.$

\medskip

\emph{First limit ordinal.} Let us denote by $\omega$ the first infinite
ordinal and by $\omega^{+}\equiv\omega+1$ its successor. We first consider the
case where
\begin{equation}
\underset{n\rightarrow\infty}{\lim\inf}\,T[g](x_{n})\geq\delta
>0.\label{eq:kra0}
\end{equation}
It follows from (\ref{eq:fat-boy}) that the sequence $\{x_{n}\}_{n\geq1}$ is
Cauchy (see Remark~\ref{rem-critical}), therefore it converges to some point $\bar{x}\in X.$ Setting
$x_{\omega}:=\bar{x}$ we deduce easily from (\ref{eq:boris1}) that
\begin{equation}
g(x_{n})>g(\bar{x}),\qquad\text{for every }n\geq0.\label{eq:kra}
\end{equation}
Similarly, we deduce from (\ref{eq:boris0}) that
\begin{equation}
\left(  f-(1+\rho)g\right)  (\bar{x})>0,\label{eq:kra1}
\end{equation}
therefore, $x_{\omega}\in X\setminus\mathcal{Z}_{T}(g)$, $\delta(x_{\omega
})>0$ and the basic iteration scheme can be pursued from $x_{\omega}$ to
define $x_{\omega+1},$ $x_{\omega+2}$ etc.\smallskip\newline We now focus on
the case
\begin{equation}
\underset{n\rightarrow\infty}{\lim\inf}\,T[g](x_{n})=0.\label{eq:kra00}
\end{equation}
Then, $x_{\omega}$ should be selected inside the set of accumulation points of
the sequence $\{x_{n}\}_{n\geq0}$ (as we have seen, this set is nonempty, but
if (\ref{eq:kra00}) holds, it might not be a singleton). Although any
accumulation point $\bar{x}$ satisfies inequalities (\ref{eq:kra}) and
(\ref{eq:kra1}), we cannot assign $x_{\omega}$ randomly among the accumulation
points, but instead, we need to select it in an adequate way (for
reasons that relate to forthcoming property (P3)  required in our forthcoming transfinite induction). To this end,
let us set $k_{0}=0$ and define inductively:
\begin{equation}
k_{n+1}:=\min\,\{m\geq k_{n}:\;T[g](x_{m})<T[g](x_{k_{n}}
)\}.\label{eq:vilches1}
\end{equation}
With this construction, $\{T[g](x_{k_n})\}_{n\in\N}$ is strictly decreasing and $T[g](x_{k_n})\to 0$. Since for all $\ell\in\lbrack k_{n},k_{n+1})\cap\mathbb{N}$ we have
$T[g](x_{\ell})\geq T[g](x_{k_{n}})$ we deduce easily from the triangular
inequality that:
\[
T[g](x_{k_{n}})\,d(x_{k_{n}},x_{k_{n+1}})\,\leq\,
{\displaystyle\sum\limits_{\ell=k_{n}}^{k_{n+1}-1}}
\,T[g](x_{\ell})\,d(x_{\ell},x_{\ell+1})
\]
yielding that
\begin{equation}
{\displaystyle\sum\limits_{n=0}^{+\infty}}
T[g](x_{k_{n}})\,d(x_{k_{n}},x_{k_{n+1}})\,\leq\,\sigma(\omega):=
{\displaystyle\sum\limits_{n=0}^{\infty}}
T[g](x_{n})\,d(x_{n},x_{n+1})\,<\,+\infty.\label{eq:vilches3}
\end{equation}
Therefore, the obtained subsequence $\{x_{k_{n}}\}_{n\geq1}$ should also have
accumulation points (it cannot be $T[g]$--critical, thanks to (iii) and
(\ref{eq:boris0})) and we define $x_{\omega}$ to be any accumulation point
$\bar{x}$ of $\{x_{k_{n}}\}_{n\geq1}$.\bigskip

\emph{Idea of the proof. }Let us outline the main idea of the proof: so far we
have defined $\{x_{n}\}_{n<\omega^{+}}\equiv\{x_{n}\}_{n<\omega}
\cup\{x_{\omega}\}$ in $X\setminus\mathcal{Z}_{T}(g)$ such that $\{g(x_{\alpha
})\}_{\alpha<\omega^{+}}$ is strictly decreasing (in terms of ordinals). Let
us denote by $\Omega$ the first uncountable ordinal, that is,
\[
\Omega=\{\,\lambda:\;\,\lambda\;\text{countable ordinal\thinspace}\}.
\]
Our objective is to extend $\{x_{n}\}_{n<\omega^{+}}$ to all countable
ordinals and come up with a generalized sequence $\{x_{\lambda}\}_{\lambda
<\Omega}$ in $X\setminus\mathcal{Z}_{T}(g)$ such that $\{g(x_{\lambda
})\}_{\lambda<\Omega}$ is strictly decreasing. Then for every $\lambda<\Omega$
we would have $g(x_{\lambda})-g(x_{\lambda^{+}})>0$ (where $\lambda^{+}$
denotes the successor of $\lambda$) and since $\Omega$ is uncountable we would
deduce:
\begin{equation}
g(x_{0})-\inf g\geq g(x_{0})-\inf_{\lambda<\Omega}g(x_{\lambda})\geq
\sum_{\lambda<\Omega}\left[  g(x_{\lambda})-g(x_{\lambda^{+}})\right]
=+\infty.\label{eq:contra}
\end{equation}
The above clearly contradicts the fact that the function $g$ is bounded and
proves the result. This construction is naturally based on transfinite
induction, where we should (also) deal with ordinals of limit type (ie.
$\lambda=\sup\{\alpha:\alpha<\lambda\}$). In this case, $x_{\lambda}$ should
be defined among the accumulation points of $\{x_{\alpha}\}_{\alpha<\lambda}$
so that we can deduce $x_{\lambda}\notin\mathcal{Z}_{T}(g)$ and guarantee the
strict descent of $g$. However, in absence of compactness, we need an
additional argument to ensure that the set accumulation points is nonempty
whenever
\[
\underset{\alpha<\lambda}{\lim\inf}\,T[g](x_{\alpha})=0\qquad\text{(else, the
limit }\underset{\alpha\nearrow\lambda}{\lim}x_{\alpha}\text{ exists!)}
\]
In this critical situation, we need to construct a sequence $\left\{
x_{\alpha_{n}}\right\}  _{n\geq1}$ with $\alpha_{n}\nearrow\lambda$ (out of
the generalized sequence $\{x_{\alpha}\}_{\alpha<\lambda}$) satisfying
\[
{\displaystyle\sum\limits_{n=0}^{+\infty}}
T[g](x_{\alpha_{n}})\,d(x_{\alpha_{n}},x_{\alpha_{n+1}})<+\infty,
\]
use assumption (iii) to deduce that $\left\{  x_{\alpha_{n}}\right\}
_{n\geq1}$ cannot be $T[g]$--critical, evoke Definition~\ref{def:AC-1} to
conclude that it has accumulation points and eventually select $x_{\lambda}$
among them. (Notice that although the set $\Omega$ is uncountable ---which is
crucial for the contradiction in (\ref{eq:contra}) above--- all of its
elements $\lambda<\Omega$ are countable ordinals and consequently, we can
always obtain cofinal sequences.) This being said, in order to effectively
realize the aforementioned critical step, we need to verify properties
(P1)--(P3) below at every step of the forthcoming transfinite induction.
\bigskip

\emph{Properties} (P1)--(P3): For every ordinal $\lambda\in\lbrack
\omega,\Omega)$, we are going to construct a generalized sequence $\{x_{\alpha
}\}_{\alpha<\lambda}\ $satisfying the following properties:

\begin{itemize}
\item[(P1)] For every $0\leq\alpha<\alpha^{+}<\lambda:$
\begin{equation}
0<\delta(x_{\alpha})d(x_{\alpha},x_{\alpha^{+}})<g(x_{\alpha})-g(x_{\alpha
^{+}}). \label{eq:fat1}
\end{equation}

\end{itemize}

Notice that (\ref{eq:fat1}) yields in particular that
\begin{equation}
\sigma(\lambda):=\,
{\displaystyle\sum\limits_{0\leq\alpha<\lambda}}
T[g](x_{\alpha})\,d(x_{\alpha},x_{\alpha^{+}})\,<\,+\infty. \label{eq:fat2a}
\end{equation}

\begin{itemize}
\item[(P2)] For every $0\leq\alpha_{1}<\alpha_{2}<\lambda:$
\begin{equation}
0\leq\left(  f-(1+\rho)g\right)  (x_{\alpha_{1}})\leq\left(  f-(1+\rho
)g\right)  (x_{\alpha_{2}}) \label{eq:fat3b}
\end{equation}
and
\begin{equation}
d(x_{\alpha_{1}},x_{\alpha_{2}})\leq
{\displaystyle\sum\limits_{\alpha_{1}\leq\alpha<\alpha_{2}}}
d(x_{\alpha},x_{\alpha^{+}}). \label{eq:fat3a}
\end{equation}

\item[(P3)] For every $\varepsilon>0$ and ordinals $0\leq\alpha_{0}
<\xi<\lambda$, there exists a finite sequence of ordinals $\alpha_{0}
<\alpha_{1}<\ldots<\alpha_{N}<\alpha_{N+1}:=\xi$ such that
\begin{equation}
{\displaystyle\sum\limits_{n=0}^{N}}
T[g](x_{\alpha_{n}})\,d(x_{\alpha_{n}},x_{\alpha_{n+1}})<\sigma(\xi
)-\sigma(\alpha_{0})+\varepsilon.\label{eq:fat4}
\end{equation}

\end{itemize}

\smallskip

\textbf{Construction (via transfinite induction).} We use transfinite
induction as follows: assuming $\{x_{\alpha}\}_{\alpha<\lambda}\subset
X\setminus\mathcal{Z}_{T}(g)$ is well-defined and satisfies (P1)--(P3), we
define $x_{\lambda}\in X\setminus\mathcal{Z}_{T}(g)$ in a way that the
extended generalized sequence $\{x_{\alpha}\}_{\alpha<\lambda^{+}}
\equiv\{x_{\alpha}\}_{\alpha\leq\lambda}$ still satisfies the same properties.
In case of a successor ordinal $\lambda=\beta^{+}$, since $x_{\beta}
\notin\mathcal{Z}_{T}(g),$ defining $x_{\lambda}\equiv x_{\beta^{+}}$ by means
of $(C)$ of Definition~\ref{def:Tmetric} automatically guarantees that
(P1)--(P3) continue to hold for $\{x_{\alpha}\}_{\alpha<\lambda^{+}}$ (see
details below). If $\lambda$ is a limit-ordinal and $x_{\lambda}$ is an
accumulation point of the generalized sequence $\{x_{\alpha}\}_{\alpha
<\lambda},$ that is,
\begin{equation}
x_{\lambda}\in
{\displaystyle\bigcap\limits_{\alpha<\lambda}}
\overline{\left\{  x_{\alpha^{\prime}}:\,\alpha^{\prime}\geq\alpha\right\}
}\qquad\text{(equivalently, }x_{\lambda}=\underset{n\rightarrow\infty}{\lim
}x_{\alpha_{n}},\;\text{for some }\alpha_{n}\nearrow\lambda\text{)}
\label{eq:fat2b}
\end{equation}
then (P1) is automatically fulfilled by the induction step (since no new
succesor ordinal is added) and (P2) follows by passing to the limit.
Therefore, the main technical difficulty is to show that $\{x_{\alpha}\}_{\alpha<\lambda}$ 
has accumulation points and to define $x_{\lambda}$
(among these accumulation points) in a way that (P3) holds. Let us now proceed to
a rigorous construction:\medskip\newline\textbf{Initialization.} We have
already defined $\{x_{\alpha}\}_{\alpha<\omega^{+}} = \{x_n\}_{n\geq 0}$ satisfying (P1)--(P2).
Let us prove that (P3) also holds. Fix $\varepsilon>0$ and consider the case
$\alpha_{0}=0$ and $\xi=\omega$ (every other value of $\alpha_0\in\N$ follows replacing $\{x_n\}_{n\geq 0}$ by $\{x_n\}_{n\geq \alpha_0}$). We consider three cases:\medskip\newline
\textit{Case} 1: $\underset{n\rightarrow\infty}{\lim\inf}\,T[g](x_{n}
)=\delta\in(0,+\infty).$ Since (\ref{eq:fat-boy}) holds, in view of Remark~\ref{rem-critical}, the sequence
$\{x_{n}\}_{n\geq1}$ converges and $x_{\omega}=\,\underset{n\rightarrow
\omega}{\lim}x_{n}.$ Taking $n_{0}\geq 0$ such that
\[
d(x_{n},x_{\omega})<\frac{\varepsilon}{\delta+1},\qquad\text{for all }n\geq
n_{0},
\]
and choosing $N\geq n_{0}$ such that $T[g](x_{N})<\delta+1$, we readily obtain
\[
{\displaystyle}\underbrace{\sum_{n=0}^{N-1}T[g](x_{n})\,d(x_{n},x_{n+1}
)}_{<\sigma(\omega)}\,+\,T[g](x_{N})\,d(x_{N},x_{\omega})\,<\,\sigma
(\omega)\,+\,\varepsilon.
\]
The result follows for $\alpha_{n}:=n,$ $n\in\{0,1,\ldots,N\}$ and
$\alpha_{N+1}:=\omega$.\medskip\newline\textit{Case} 2:
$\underset{n\rightarrow\infty}{\lim\inf}\,T[g](x_{n})=+\infty.$ Similarly to
the previous case, the sequence $\{x_{n}\}_{n\geq1}$ converges and $x_{\omega
}=\,\underset{n\rightarrow\omega}{\lim}x_{n}.$ Let $N\geq0$ be such that
\[
T[g](x_{N})\,=\,\underset{n\geq0}{\min}\{T[g](x_{n})\}\quad\text{(the minimum
is attained since }T[g](x_{n})\rightarrow+\infty\text{).}
\]
Then, $T[g](x_{n})\geq T[g](x_{N})$ for all $n\geq N$ which in view of
(\ref{eq:fat3a}) (for $\alpha_{1}=N$ and $\alpha_{2}=\omega$) yields
\[
T[g](x_{N})\,d(x_{N},x_{\omega})\,\leq\,T[g](x_{N})\,
{\displaystyle\sum\limits_{n=N}^{+\infty}}
d(x_{n},x_{n+1})<\sum_{n=N}^{+\infty}T[g](x_{n})\,d(x_{n},x_{n+1}).
\]
Consequently,
\[
\sum_{n=0}^{N-1}T[g](x_{n})\,d(x_{n},x_{n+1})\,+\,T[g](x_{N})\,d(x_{N}
,x_{\omega})\,\leq\,\sum_{n=0}^{+\infty}T[g](x_{n})\,d(x_{n},x_{n+1}
):=\sigma(\omega).
\]

\textit{Case }3: $\underset{n\rightarrow\infty}{\lim\inf}\,T[g](x_{n})=0.$ In
this case $x_{\omega}$ is defined among the accumulation points of the
sequence $\{x_{k_{n}}\}_{n\geq1}$ constructed in (\ref{eq:vilches1}). Given
$\varepsilon>0$ we chose $n_{0}\in\mathbb{N}$ in a way that 
\[
d(x_{k_{n_{0}}},x_{\omega})\,<\,\frac{\varepsilon}{T[g](x_{k_{n_0}})}. 
\]
In view of \eqref{eq:vilches3}, we deduce that
\[
{\displaystyle\sum\limits_{n=0}^{n_{0}-1}}
T[g](x_{k_{n}})\,d(x_{k_{n}},x_{k_{n+1}})\,+\,T[g](x_{k_{n_{0}}
})\,d(x_{k_{n_{0}}},x_{\omega})\,<\,\sigma(\omega)\,+\,\varepsilon,
\]
and the result follows for $N=k_{n_{0}}$ and $\alpha_{n}:=k_{n},$
$n\in\{0,1,\ldots,N\}.$\medskip\newline\textbf{Successor ordinal.} Assume
$\lambda=\beta^{+}>\omega$ is a successor ordinal and $\{x_{\alpha}
\}_{\alpha<\lambda}\equiv\{x_{\alpha}\}_{\alpha\leq\beta}$ is well-defined and
satisfies (P1)--(P3). Then from (\ref{eq:fat3b}) and (ii) we deduce that
$x_{\beta}\notin\mathcal{Z}_{T}(g)$ and $\delta(x_{\beta})>0.$ Using $(C)$ as
before, we obtain $x_{\beta^{+}}\in\lbrack g<g(x_{\beta})]$ such that
\[
0<f(x_{\beta})-(1+\rho)g(x_{\beta})\leq f(x_{\beta^{+}})-(1+\rho
)g(x_{\beta^{+}})\quad\text{and\quad}\delta(x_{\beta})d(x_{\beta},x_{\beta
^{+}})\leq g(x_{\beta})-g(x_{\beta^{+}}).
\]
Notice that (\ref{eq:fat3a}), (\ref{eq:fat4}) follow easily from the
induction step and the triangular inequality. Therefore $\{x_{\alpha
}\}_{\alpha<\lambda^{+}}\equiv\{x_{\alpha}\}_{\alpha\leq\lambda}$ also
satisfies (P1)--(P3).\bigskip\newline\textbf{Ordinal of limit type.} Let us
now assume that $\lambda\in(\omega,\Omega)$ is a limit ordinal and
$\{x_{\alpha}\}_{\alpha<\lambda}$ is defined and satisfies (P1)--(P3).
\smallskip

\emph{Case} \emph{I}. We first focus on the case where
\begin{equation}
\underset{\alpha<\lambda}{\lim\inf}\,T[g](x_{\alpha})\geq\delta
>0.\label{eq:amaya}
\end{equation}
We deduce that there exists $\xi_0<\lambda$ such that $T[g](x_{\alpha})\geq \delta/2$ for all $\xi_0\leq\alpha<\lambda$. Thus, in an analogous way as in Remark~\ref{rem-critical}, we get from (\ref{eq:fat2a}) that
\[
{\displaystyle\sum\limits_{\xi_0\leq\alpha<\lambda}}
\,d(x_{\alpha},x_{\alpha^{+}})\,\leq\,\frac{2}{\delta}
{\displaystyle\sum\limits_{\xi_0\leq\alpha<\lambda}}
T[g](x_{\alpha})\,d(x_{\alpha},x_{\alpha^{+}})\,<\,+\infty
\]
and consequently the generalized sequence $\{x_{\alpha}\}_{\alpha<\lambda}$
converges as $\alpha\nearrow\lambda.$ We set
\[
x_{\lambda}:=\underset{\alpha\rightarrow\lambda}{\lim}x_{\alpha}
\]
and obtain readily that the generalized sequence $\{x_{\alpha}\}_{\alpha
<\lambda^{+}}$ still satisfies (P1)--(P2). It remains to check that (P3) holds
for $\varepsilon>0$ and $0\leq\alpha_{0}<\lambda<\lambda^{+}.$ To this end, we
need to consider separately two different cases depending on whether
(\ref{eq:amaya}) is finite or not.\smallskip\newline-- \emph{Subcase}
$I_{1}$. Assume that $\underset{\alpha<\lambda}{\lim\inf}
\,T[g](x_{\alpha})=\delta\in(0,+\infty).$ Then, we fix $\alpha_{\ast}\geq
\alpha_{0}$ such that for all $\alpha\in(\alpha_{\ast},\lambda)$ we have
\[
d(x_{\alpha},x_{\lambda})<\frac{\varepsilon}{2(\delta+1)}.
\]
Pick any $\hat{\alpha}\in(\alpha_{\ast},\lambda)$ such that $T[g](x_{\hat
{\alpha}})<\delta+1$. We deduce directly that
\[
T[g](x_{\hat{\alpha}})\,d(x_{\hat{\alpha}},x_{\lambda})\,<\,\frac{\varepsilon}{2}.
\]
Applying (P3) for $0\leq\alpha_{0}<\hat{\alpha}<\lambda$ (and $\tilde
{\varepsilon}=\varepsilon/2$) we obtain a finite strictly increasing sequence
$\{\alpha_{n}\}_{n=0}^{N}$ satisfying
\[
{\displaystyle\sum\limits_{n=0}^{N-1}}
T[g](x_{\alpha_{n}})\,d(x_{\alpha_{n}},x_{\alpha_{n+1}})+T[g](x_{\alpha_{N}
})\,d(x_{\alpha_{N}},x_{\hat{\alpha}})<\sigma(\hat{\alpha})-\sigma(\alpha
_{0})+\frac{\varepsilon}{2}.
\]
The result follows by concatenation.\smallskip\newline-- \emph{Subcase}
\emph{I}$_{2}$. Assume that $\underset{\alpha<\lambda}{\lim\inf}
\,T[g](x_{\alpha})=\,\underset{\alpha<\lambda}{\lim}\,T[g](x_{\alpha}
)=+\infty.$ Then, there exists $\alpha_{\ast}\geq\alpha_{0}$ such that
$T[g](x_{\alpha})\geq 1,$ for all $\alpha\in(\alpha_{\ast},\lambda)$. We deduce from
(\ref{eq:fat2a}) that
\[
\Delta:=
{\displaystyle\sum\limits_{\alpha_{\ast}\leq\alpha<\lambda}}
d(x_{\alpha},x_{\alpha^{+}})\leq
{\displaystyle\sum\limits_{\alpha_{\ast}\leq\alpha<\lambda}}
T[g](x_{\alpha})\,d(x_{\alpha},x_{\alpha^{+}})<+\infty.
\]
Set
\[
\mu_{\ast}:=\,\underset{\alpha_{\ast}\leq\alpha<\lambda}{\inf}\{T[g](x_{\alpha
})\}\geq1.
\]
In contrast to Case 2 of Initialization, the above infimum might not be attained if $\lambda$ is limit of ordinals of limit
type (that is, $\lambda=\sup\{\xi<\lambda:\,\xi$ limit-ordinal$\}$). However,
we can choose $\hat{\alpha}\in(a_{\ast},\lambda)$ such that
\[
T[g](x_{\hat{\alpha}})<\mu_{\ast}+\frac{\varepsilon}{2\Delta}.
\]
We deduce from (\ref{eq:fat3a}) (for $\alpha_{1}=\hat{\alpha}$ and $\alpha
_{2}=\lambda$) that
\[
\begin{aligned}
T[g](x_{\hat{\alpha}})d(x_{\hat{\alpha}},x_{\lambda})\leq T[g](x_{\hat{\alpha
}})
{\displaystyle\sum\limits_{\hat{\alpha}\leq\alpha<\lambda}}
d(x_{\alpha},x_{\alpha^{+}})\,<\,&
{\displaystyle\sum\limits_{\hat{\alpha}\leq\alpha<\lambda}}
\left(  \mu_{\ast}+\frac{\varepsilon}{2\Delta}\right)  \,d(x_{\alpha
},x_{\alpha^{+}}) \\ <\,&\underbrace{
{\displaystyle\sum\limits_{\hat{\alpha}\leq\alpha<\lambda}}
T[g](x_{\alpha})\,d(x_{\alpha},x_{\alpha^{+}})}_{:=\sigma(\lambda)-\sigma
(\hat{\alpha})}\,+\,\frac{\varepsilon}{2}.
\end{aligned}
\]
We conclude as before by concatenating $\{\hat{\alpha},x_{\lambda} \}$ 
with the finite sequence obtained by
applying~(P3) (which by the induction step is assumed to hold for the generalized sequence $\{x_{\alpha}\}_{\alpha < \lambda}$) for the choice $\tilde
{\varepsilon}=\varepsilon/2$ and $0\leq\alpha_{0}<\hat{\alpha}<\lambda$. \medskip \newline
\emph{Case} \emph{II}. It remains to deal with the case
\[
\underset{\alpha<\lambda}{\lim\inf}\,T[g](x_{\alpha})=0.
\]
Our objective is to show that
\[
\mathcal{C}=
{\displaystyle\bigcap\limits_{\alpha<\lambda}}
\overline{\left\{  x_{\alpha^{\prime}}:\,\alpha^{\prime}\geq\alpha\right\}
}\neq\emptyset
\]
and define $x_{\lambda}\ $in $\mathcal{C}$ in a way that $\{x_{\alpha
}\}_{\alpha<\lambda^{+}}$ satifies (P3). (We recall that if $x_{\lambda}
\in\mathcal{C}$ then (P1)--(P2) are automatically satisfied.) To this end, fix
$\varepsilon>0$ and $0\leq\alpha_{0}<\lambda<\lambda^{+}.$ Let
\[
\{\varepsilon_{n}\}_{n\geq0}\subset(0,\varepsilon)\qquad\text{such that}\quad
{\displaystyle\sum\limits_{n=0}^{+\infty}}
\varepsilon_{n}=\frac{\varepsilon}{2}.
\]
Let further $\gamma_{n}\nearrow\lambda$ and define inductively a sequence of
ordinals $\{\xi_{n}\}_{n\geq0}$ as follows:
\[
\xi_{0}=\alpha_{0}\qquad\text{and}\qquad\xi_{n+1}:=\min\,\left\{  \alpha
\geq\max\{\xi_{n},\gamma_{n}\}:\,T[g](x_{\alpha})\leq\frac{T[g](x_{\xi_{n}}
)}{2}\right\}  .
\]
The above definition guarantees that $\xi_{n}\nearrow\lambda$ and
$T[g](x_{\xi_{n}})\searrow0$ (that is, the sequence of descent moduli
converges to $0$ decreasingly). For every $n\geq0$, thanks to our induction assumption, 
we can apply property (P3) for
$\varepsilon_{n}>0$ and $0\leq\xi_{n}<\xi_{n+1}<\lambda^{+}$ to obtain
$N_{n}\geq1$ and a finite sequence 
$\xi_{n}:=\alpha_{0}^{n}<\alpha_{1}^{n}<\ldots<\alpha_{N_{n}+1}^{n}:=\xi_{n+1}$
such that
\[
{\displaystyle\sum\limits_{i=0}^{N_{n}}}
T[g](x_{\alpha_{i}^{n}})\,d(x_{\alpha_{i}^{n}},x_{\alpha_{i+1}^{n}
})\,<\,\sigma(\xi_{n+1})-\sigma(\xi_{n})+\varepsilon_{n}.
\]
Concatenating the above finite sequences $\{{\alpha}_i^n:\,i\in\{0,\dots,N_n\}\}$, for $n\geq 0$, we obtain a strictly increasing sequence
\begin{equation}
\alpha_{0}:=\alpha_{0}^{0}<\alpha_{1}^{0}<\ldots<\alpha_{N_{0}+1}^{0}:=\xi
_{1}:=\alpha_{0}^{1}<\ldots<\alpha_{N_{1}+1}^{1}:=\xi_{2}<\ldots<\ldots
\label{eq:amaya3}
\end{equation}
which converges to $\lambda$ and satisfies
\[
{\displaystyle\sum\limits_{n=0}^{+\infty}}
{\displaystyle\sum\limits_{i=0}^{N_{n}}}
T[g](x_{\alpha_{i}^{n}})\,d(x_{\alpha_{i}^{n}},x_{\alpha_{i+1}^{n}
})\,<\,\sigma(\lambda)-\sigma(\alpha_{0})+\frac{\varepsilon}{2}.
\]
Renaming (\ref{eq:amaya3}) to $\{\beta_{n}\}_{n\geq0}$, we have $\beta
_{0}\equiv\alpha_{0},$ $\beta_{n}\nearrow\lambda$
\[
{\displaystyle\sum\limits_{n=0}^{+\infty}}
T[g](x_{\beta_{n}})\,d(x_{\beta_{n}},x_{\beta_{n}+1})\,<\,\sigma
(\lambda)-\sigma(\alpha_{0})+\frac{\varepsilon}{2}\quad\text{and}
\quad\underset{n\rightarrow\infty}{\lim\inf}\,T[g](x_{\beta_{n}})=0.
\]
Acting as in (\ref{eq:vilches1}), we set $k_{0}=\beta_{0}\equiv\alpha_{0}$ and
\[
k_{n+1}:=\min\,\{\beta_m\geq k_{n}:\;T[g](x_{\beta_{m}})<T[g](x_{{k_{n}}})\}.
\]
Since for all $\ell\in\lbrack k_{n},k_{n+1})\cap\{\beta_m\}_{m\in\mathbb{N}}$ we have
$T[g](x_{\beta_{\ell}})\geq T[g](x_{\beta_{k_{n}}})$, we deduce:
\begin{equation}
{\displaystyle\sum\limits_{n=0}^{+\infty}}
T[g](x_{k_{n}})\,d(x_{k_{n}},x_{k_{n+1}})\,\leq\,
{\displaystyle\sum\limits_{n=0}^{+\infty}}
{\displaystyle\sum\limits_{\ell=k_{n}}^{k_{n+1}-1}}
\,T[g](x_{\ell})\,d(x_{\ell},x_{\ell+1})<\,\sigma(\lambda)-\sigma(\alpha
_{0})+\frac{\varepsilon}{2}.\label{salas1}
\end{equation}
Since $\{x_{k_{n}}\}_{n\geq1}$ cannot be $T[g]$--critical (thanks to
assumption (iii)), we deduce from Definition~\ref{def:Tmetric} that it has
accumulation points as $n\rightarrow\infty.$ We define $x_{\lambda}$ to be any
accumulation point of $\{x_{k_{n}}\}_{n\geq1}.$ Then, the last part of the
argument is the same as in Case 3 of Initialization: we chose $n_{0}\in\mathbb{N}$ in a way that
\[
d(x_{k_{n_{0}}},x_{\lambda})\,<\,\frac{\varepsilon}{2T[g](x_{k_{0}})}\,
\]
$\ $and we deduce that
\begin{equation}
{\displaystyle\sum\limits_{n=0}^{n_{0}-1}}
T[g](x_{k_{n}})\,d(x_{k_{n}},x_{ k_{n+1}})\,+\,T[g](x_{k_{n_{0}}
})\,d(x_{k_{n_{0}}},x_{\lambda})\,<\,
{\displaystyle\sum\limits_{n=0}^{+\infty}}
T[g](x_{k_{n}})\,d(x_{k_{n}},x_{k_{n+1}})\,+\,\frac{\varepsilon}
{2},\label{salas2}
\end{equation}
and the result follows by combining (\ref{salas1}) and (\ref{salas2}).\bigskip\newline This completes the transfinite induction and allows to
obtain a generalized sequence $\{x_{\alpha}\}_{\alpha<\lambda}$ for all
countable ordinal $\lambda$. Thus we ultimately define $\{x_{\alpha
}\}_{\alpha<\Omega}$ with $\{g(x_{\lambda})\}_{\lambda<\Omega}$ uncountable
and strictly decreasing, contradicting that $g$ is bounded from
below.\smallskip\newline Therefore, $f<(1+\rho)g.$ Repeating the procedure for
any $\rho>0$, we deduce that $f\leq g$ should hold. The proof is
complete.\hfill$\square$

\bigskip

\begin{remark}
[discussion on the assumptions]\label{rem-orteg}$(a)$. Assumptions (ii) and
(iii) are complementary and independent: indeed, asymptotically critical
sequences can neither yield nor be obtained by critical points, since they are
not allowed to converge (c.f. Definition~\ref{def:AC-1}). In particular, in
absence of compactness, the set of critical points $\mathcal{Z}_{T}(g)$ can be
empty, case in which assumption (ii) of Theorem~\ref{thm.strict_compa} is
trivially satisfied and provides no information. This potential lack of
information is contemplated by assumption~(iii).\smallskip\newline
$(b)$. The main difficulty in proving Theorem~\ref{thm.strict_compa} is that the construction requires transfinite induction, while assumption~(iii) allows comparison through sequences.  Indeed, from a point $x_{\alpha}$ for which $T[g](x_{\alpha})>0$, we produce a descent point $x_{\alpha^+}$ (that is, $g(x_{\alpha^+}) < g(x_{\alpha})$, and after countably many descent points, we select an adequate accumulation point; due to the (possible) existence of points $x\in X\setminus \mathcal{Z}_{T}(g)$ for which $\liminf_{y\to x} T[g](y) = 0$, such construction might end prematurely unless we allow to restart the process from the limit points, inducing a transfinite construction (see Example~\ref{ex:ManyLowerCrit}). The critical part is to prove, using only asymptotically critical sequences, that for every limiting ordinal $\lambda$, the accumulation point $x_{\lambda}$ can always be constructed. The importance of the invariant properties (P1)--(P3) during this construction was precisely the fact that if such a point $x_{\lambda}$ fails to exist, then we would be able to extract an asymptotically critical (cofinal) subsequence from $\{x_{\alpha}\}_{\alpha<\lambda}$, and compare the functions over that sequence.
\smallskip\newline
$(c)$. Let
us momentarily assume that condition (\ref{eq:morduk}) is imposed to every
sequence $\{z_{n}\}_{n\geq1}\subset X\setminus\mathcal{Z}_{T}(g)$ satisfying
(\ref{eq:AC}), rather than only to those that are free of accumulation points.
Let further $\bar{z}$ denote some accumulation point of $\{z_{n}\}_{n\geq1}.$
Then the case where $\bar{z}$ is critical (i.e. $\bar{z}\in\mathcal{Z}_T(g)$) is
already covered by assumption~(ii) of Theorem~\ref{thm.strict_compa} (since
$f$ is continuous and $g$ lower semicontinuous) while the case where $\bar
{z}\notin\mathcal{Z}_T(g)$ leads to a superfluous assumption making the
statement of the theorem weaker. This is illustrated in Example
\ref{ex:ManyLowerCrit} below, where the descent operator $T[f]$ is the local
slope $s[f]$. The example reveals that accumulation points of sequences
satisfying (\ref{eq:AC}) might not be critical for the slope $s[g]$ but
instead, for the closure of $s[g]$ (called regularized slope in~\cite{DLI2015}) and that it is neither necessary nor desirable to impose any condition there.
\end{remark}

\begin{example}
\label{ex:ManyLowerCrit} Set $X=[1,+\infty)$ with the usual distance. For each
interval $I_{n}=[n,n+1)$ we define the function $g_{n}:I_{n}\rightarrow
[0,1]$ given by
\[
g_{n}(x)=\frac{1}{n+1}+\left(  \frac{1}{n(n+1)}\right)  \left(
x-(n+1)\right)^{n(n+1)}.
\]
We finally define $g:[1,+\infty)\rightarrow [0,1]$ given by $g(x)=g_{n}
(x)$ whenever $x\in I_{n}$.
\begin{figure}[!h]
    \centering
    \begin{tikzpicture}[scale=0.8]
    \begin{axis}[
        grid = major,
        grid style={dashed, gray!50},
        scale only axis,
        xmin=1,xmax=11,
        xtick={1,2,4,6,8,10},
    ]
    {\foreach \nn in {1,2,...,9}
        \addplot[black,thick, domain = \nn:\nn+1]{1/(\nn+1) + (1/(\nn*(\nn+1)))*(x-\nn-1)^(\nn*(\nn+1)) };
    }
    \addplot[dotted,black,thick,domain=10:11]{1/(11) + (1/(11*12))*(x-10-1)^2 };
    \addplot [red,thick, mark=o, only marks] coordinates { (1, {1/2 + 1/(1*2)} )(2, {1/3 + 1/(2*3)} ) 
    (3, {1/4 + 1/(3*4)} )(4, {1/5 + 1/(4*5)} )
    (5, {1/6 + 1/(5*6)} )(6, {1/7 + 1/(6*7)} ) 
    (7, {1/8 + 1/(7*8)} )(8, {1/9 + 1/(8*9)} )
    (9, {1/10 + 1/(9*10)} )(10, {1/11 + 1/(10*11)} )
    };
    \end{axis}
    \end{tikzpicture}
    \caption{ Function $g:[1,+\infty)\to [0,1]$, constructed by blocks $I_n=[n,n+1)$. Plot of the first 9 intervals. At each integer point $n\geq 1$, the lateral derivatives are $g^{\prime}_{-}(n)=0$ from the left and $g^{\prime}_{+}(n)=1$ from the right. }\label{fig:LowerCritical}
\end{figure}
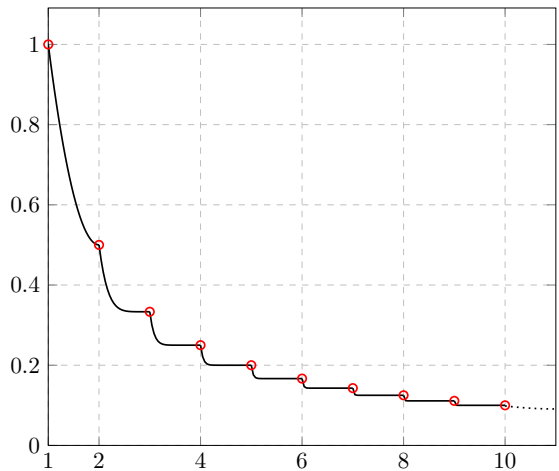

Consider the descent operator $T[g]=s[g]$ (local
slope) and notice that
\[
s[g](x)=(n+1-x)^{n(n+1)-1}>0,\quad\text{for }x\in\lbrack n,n+1)\qquad\text{
and }\qquad\mathcal{Z}_{T}[g]=\emptyset.
\]
Notice further that
\[
\underset{x\in X}{\inf}g=\lim_{x\rightarrow+\infty}g(x)=0.
\]
Finally, for every $n\in\mathbb{N}$, the point $\bar{x}:=n$ is not critical,
but it is critical for the regularized slope (see Figure~\ref{fig:LowerCritical}), that is:
\begin{equation}
\underset{x\rightarrow n}{\liminf}\,s[g](x)=0. \label{eq:Boris}
\end{equation}
\smallskip\newline
It is easy to see that $f\leq g$ for every continuous function $f:X\rightarrow
\mathbb{R}$ satisfying
\[
(a)\quad s[f]<s[g]\,\,\,\text{ on }\,X\qquad\text{ and }\qquad(b)\quad
\liminf_{x\rightarrow+\infty}f(x)<0.
\]
Indeed, we can either apply Theorem~\ref{thm.strict_compa} or do the following
elementary proof: pick any increasing sequence $x_{n}\rightarrow+\infty$ such
that $\underset{n\rightarrow+\infty}{\lim}f(x_{n})=\underset{x\rightarrow
+\infty}{\liminf}f(x)$, and assume towards a contradiction that $f(x_{0}%
)\geq(1+\rho)g(x_{0})$, for some $\rho>0$. Following the construction of
\cite[Lemma~3.3]{DMS2022}, we can build a (generalized) sequence
$\{z_{\lambda}\}_{\lambda}$ in the compact interval $[x_{0},x_{1}]$, strictly
decreasing for $g$ and such that $f(z_{\lambda})\geq(1+\rho)g(z_{\lambda})$
for each $\lambda$. The construction eventually ends (due to cardinality
obstructions) and the only way for this to happen is that $z_{\lambda
}\rightarrow x_{1}$, since this is the only critical point of $g$ restricted
to $[x_{0},x_{1}]$. Continuity of $f$ and lower semicontinuity of $g$ would
yield that $f(x_{1})\geq(1+\rho)g(x_{1})$. An inductive argument shows that
$f(x_{n})\geq(1+\rho)g(x_{n})$ for all $n\in\mathbb{N}$, which yields
$\lim_{n}f(x_{n})\geq\inf g=0$, leading to a contradiction.\smallskip\newline
On the other hand, if Definition~\ref{def:AC-1} allowed to consider convergent
sequences, then Theorem~\ref{thm.strict_compa} could not directly apply since
it would have required an extra condition on all sequences satisfying
(\ref{eq:Boris}), leading to (infinitely many) unnecessary extra conditions:
$f(n)<g(n)$ for all $n\geq1$.

Finally, observe that for every $n\in\N$, the steepest descent curve $\gamma_n:[0,+\infty)\to\R$ solving 
\begin{equation*}
\begin{cases}
\dot{\gamma}(t)=-\nabla g(\gamma(t)),\quad t\geq0,\\
\gamma(0)=n,
\end{cases}
\end{equation*}
satisfies that $\displaystyle\lim_{t\to+\infty} \gamma_n(t) = n+1$. Thus, if we follow the construction of Theorem~\ref{thm.strict_compa} by taking a descent point at each iteration, we should obtain a generalized sequence $\{x_{\alpha}\}$ (similar to a concatenation of discretizations of the curves $\{\gamma_n\}$) diverging to $+\infty$. The delicate construction of the proof of Theorem~\ref{thm.strict_compa} would allow us to retrieve an $s$-asymptotically critical cofinal subsequence for $g$, and therefore to compare $f$ and $g$ at the limit values. \hfill$\Diamond$
\end{example}

\subsection{Determination in complete metric spaces}

\label{sec:MainResult}If two functions $f,g$ have the same descent modulus at
every point (that is, $T[f]=T[g]$), then they have the same critical set
($\mathcal{Z}_{T}:=\mathcal{Z}_{T}(f)=\mathcal{Z}_{T}(g)$) and the same
asymptotically critical sequences ($\mathcal{AZ}_{T}:=\mathcal{AZ}
_{T}(f)=\mathcal{AZ}_{T}(g)$). Using the same strategy as in
Theorem~\ref{thm:DetExtended-compact}, we obtain the main result of this work.

\begin{theorem}
[main determination result]\label{thm:DetExtended-metric}Let $T$ be a
metrically compatible descent modulus on the complete metric space $(X,d)$.
Let $f,g\in\overline{\mathcal{C}}(X)\cap\mathrm{dom}(T)$ be bounded from below
and satisfy
\[
T[f](x)=T[g](x),\qquad\text{for all }x\in X.
\]
Assume that
\thinspace\thinspace\ $f|_{\mathcal{Z}_{T}}=g|_{\mathcal{Z}_{T}}$\quad and
\quad$\underset{n\rightarrow\infty}{\lim\inf}\,f(z_{n}
)=\,\underset{n\rightarrow\infty}{\lim\inf}\,g(z_{n})$,\thinspace
\thinspace\ for all $\{z_{n}\}_{n\geq1}\in\mathcal{AZ}_{T}$.\medskip\newline
Then
\[
f(x)=g(x),\quad\text{for all }x\in X.
\]

\end{theorem}

By considering only metric spaces and metrically compatible descent moduli, it
is clear that the above result generalizes
Theorem~\ref{thm:DetExtended-compact}: indeed, for every $f\in\mathcal{F}$ one
has that $\mathcal{AZ}_{T}(f) = \emptyset$ whenever $X$ is compact.
\smallskip\newline

Theorem~\ref{thm:DetExtended-metric} also generalizes \cite[Section~4]
{TZ2022}. This is due to the fact that for the global slope $\mathcal{G}$
defined in~\eqref{eq:globalSlope}, every $\mathcal{G}$-asymptotically
critical sequence of $f$ 
is infimizing for the function $f$. This is
the content of the following lemma.

\begin{lemma}
[infimizing sequences]\label{lem:jortega} Let $(X,d)$ be a complete metric
space, $f\in\mathcal{\bar{C}}(X)\cap\dom(\mathcal{G})$ and $\{z_{n}\}_{n\geq
1}$ a $\mathcal{G}$-asymptotically critical sequence
(Definition~\ref{def:AC-1}). Then
\begin{equation}
\liminf_{n\rightarrow\infty}\,f(z_{n})=\inf\,f. \label{eq:criticalseq-inf}
\end{equation}

\end{lemma}

\noindent\textbf{Proof.} Let $\{z_{n}\}_{n\geq1}$ be a $\mathcal{G}
$-asymptotically critical sequence for the function $f.$ Then, $\{\mathcal{G}
[f](z_{n})\}_{n\geq1}$ is a sequence of strictly positive numbers that
converges to zero. We set $k_{1}=1$ and define inductively
\[
k_{n+1}:=\min\,\left\{  m\geq k_{n}:\;\mathcal{G}[f](z_{m})<\mathcal{G}
[f](z_{k_{n}})\right\}  .
\]
Then, $\{\mathcal{G}[f](z_{k_{n}})\}_{n\geq1}$ is strictly decreasing and for
$m\in\lbrack k_{n},k_{n+1})\cap\mathbb{N}$ we have
\[
\mathcal{G}[f](z_{k_{n}})\leq\mathcal{G}[f](z_{m}).
\]
We deduce that for every $n\geq1$
\begin{align*}
\mathcal{G}[f](z_{k_{n}})\,d(z_{k_{n}},z_{k_{n+1}})\,  &  \leq\,\mathcal{G}
[f](z_{k_{n}})\,\sum_{k_{n}\leq m<k_{n+1}}d(z_{m},z_{m+1})\\
&  \leq\sum_{k_{n}\leq m<k_{n+1}}\,\mathcal{G}[f](z_{m})d(z_{m},z_{m+1})
\end{align*}
and consequently
\[
\sum_{n=1}^{\infty}\mathcal{G}[f](z_{k_{n}})\,d(z_{k_{n}},z_{k_{n+1}}
)\,\leq\,\sum_{m=1}^{\infty}\,\mathcal{G}[f](z_{m})d(z_{m},z_{m+1})<+\infty.
\]
Let $u\in X$ be arbitrarily chosen. We deduce from the definition of the
global slope (\ref{eq:globalSlope}) that
\[
f(z_{k_{n}})\leq f(u)+\mathcal{G}[f](z_{k_{n}})\,d(z_{k_{n}},u),\quad\text{for
all }n\geq1.
\]
Therefore, it suffices to show that
\[
\liminf_{n\rightarrow\infty}\,\mathcal{G}[f](z_{k_{n}})\,d(z_{k_{n}},u)=0.
\]
To this end, take $n,m\in\mathbb{N}$ with $n<m$. Then, since $\mathcal{G}
[f](z_{k_{n}})$ is decreasing, we deduce:
\[
\mathcal{G}[f](z_{k_{m}})\,d(z_{k_{m}},u)\leq\mathcal{G}[f](z_{k_{m-1}
})\,d(z_{k_{m-1}},z_{k_{m}})+\mathcal{G}[f](z_{k_{m}})\,d(z_{k_{m-1}},u)
\]
and consequently
\[
\mathcal{G}[f](z_{k_{m}})\,d(z_{k_{m}},u)\,\leq\,\sum_{n\leq\ell<m}
\mathcal{G}[f](z_{k_{\ell}})\,d(z_{k_{\ell}},z_{k_{\ell+1}})+\mathcal{G}
[f](z_{k_{m}})\,d(z_{k_{n}},u).
\]
Thus, keeping $n$ fixed and passing to the limit as $m\rightarrow\infty$, we
deduce that:
\[
\liminf_{m\rightarrow\infty}\,\mathcal{G}[f](z_{k_{m}})\,d(z_{k_{m}}
,u)\,\leq\,\sum_{\ell=n}^{+\infty}\mathcal{G}[f](z_{k_{\ell}})\,d(z_{k_{\ell}
},z_{k_{\ell+1}})\leq\sum_{\ell=n}^{\infty}\mathcal{G}[f](z_{\ell})\,d(z_{\ell}
,z_{\ell+1}).
\]
Since the last quantity becomes arbitrarily small as $n$ increases, the
conclusion follows. \hfill$\Box$

\bigskip
\begin{remark}\label{rem-penot}
Comparing the above lemma with \cite[Lemma 4.2]{TZ2022}, we observe that the main difference is that the latter considers  summable sequences where $\sum_{n} \mathcal{G}[g](z_{n+1})d(z_n,z_{n+1}) < +\infty$. This can be seen as proximal algorithm-type condition. However, in our context, we need to consider asymptotically critical sequences where $\sum_{n} \mathcal{G}[g](z_{n})d(z_n,z_{n+1}) < +\infty$. Thus, asymptotically critical sequences can be seen as gradient algorithm-type condition, which is strongly related on how descent (generalized) sequences are constructed from descent moduli.
	\end{remark}

\smallskip
Finally, Lemma~\ref{lem:jortega}, together with the fact that any critical point for the global
slope has to be a global minimizer, yields directly the following corollary.

\begin{corollary}
[Global slope determination, {\cite{TZ2022}}]\label{cor:jortega}Let $(X,d)$ be
a complete metric space and $f,g:X\rightarrow\mathbb{R}$ be two proper lower
semicontinuous functions which are bounded from below and continuous on their
domain. If $\mathcal{G}[f](x)=\mathcal{G}[g](x)$, for all $x\in X$ and $\inf
f=\inf g$, then $f=g$.
\end{corollary}

\bigskip

Corollary~\ref{cor_ort} below reveals that a continuous bounded from below
function $f\in\mathcal{\bar{C}}(X)\cap\mathrm{dom}(T)$ necessarily possesses
either a critical point or an asymptotically critical sequence.

\begin{corollary}
[existence of critical elements]\label{cor_ort} Let $f:X\rightarrow
\mathbb{R}\cup\{+\infty\}$ be a lower semicontinuous function which is
continuous on its domain and bounded from below. If $f\in\mathrm{dom}(T)$ for
some metrically compatible descent modulus on $X$, then either $\mathcal{Z}
_{T}(f)\neq\emptyset$ or $\mathcal{AZ}_{T}(f)\neq\emptyset$.
\end{corollary}

\noindent\textbf{Proof. }We may assume that $\inf f=0.$ Fix $\varepsilon>0$,
set $g:=(1+\varepsilon)\,f$ and $\tilde{f}=f+1.$ It follows that
$\inf\,g=0\,<$ $\inf\,\tilde{f}=1.$ In view of
Proposition~\ref{prop.mono_proper} and property $(\mathcal{D}_{3})$ of the
descent modulus, we deduce that $\mathcal{Z}_{T}(g)\subset\mathcal{Z}
_{T}(\tilde{f})$ and $T[g](x)>T[\tilde{f}](x),$ for all $x\in X\setminus
\mathcal{Z}_{T}(g).$ Moreover, if a sequence $\{z_{n}\}_{n}$ is $T[g]$
--critical, then it is also $T[\tilde{f}]$--critical (and consequently,
$T[f]$--critical). Let us assume, towards a contradiction that $g$ has no
critical points and no $T[g]$--critical sequences. Then $\mathcal{Z}
_{T}(g)=\emptyset$ and assumptions (ii) and (iii) of
Theorem~\ref{thm.strict_compa} are trivially fulfilled. We deduce that
$g>\tilde{f}$ which is a contradiction. Therefore, either $\mathcal{Z}
_{T}(g)\neq\emptyset$ or $\mathcal{AZ}_{T}(g)\neq\emptyset$.\hfill$\square$

\bigskip

\noindent(\textit{Open question}) It is well-known that the local and the
global slopes coincide for convex functions in any Banach space and are equal
to the remoteness of the subdifferential (the distance of the convex
subdifferential to zero). This fact was used in \cite[Section 5]{TZ2022} to
obtain a nontrivial generalization (from Hilbert to Banach spaces) of the
determination result for the class of convex functions obtained in
\cite[Theorem 3.8]{BCD2018} (smooth case) and \cite{PSV2021} (nonsmooth case).
Indeed, the identification of derivatives with gradients and the use of
(sub)gradient systems played a crucial role in the latter works. In the recent
work \cite{DD2023}, the authors have again used the Hilbertian structure to
show that the deviation between the slopes of two convex functions controls
the deviation between the functions themselves. It is not known if such
slope-based sensitivity result would hold for convex functions in Banach
spaces, or more generally, if one can use metric descent modulus deviations to
measure deviations of functions in general.

\bigskip

\noindent\rule{5cm}{1pt} \smallskip\newline\noindent\textbf{Acknowledgements.}
The authors wish to thank Sebastian Tapia-García for useful discussions
related to Remark~\ref{rem:tap}. The first author acknowledges support from
the Austrian Science Fund (\textsc{FWF, P-36344-N}). The third author was
partially funded by \textsc{ANID}-Chile, through the grant \textsc{Fondecyt}
\textsc{11220586}, and \textsc{FB210005 BASAL} funds for centers of
excellence.

\vspace{0.5cm}

\noindent Aris DANIILIDIS, Tri Minh LE

\medskip

\noindent Institut f\"{u}r Stochastik und Wirtschaftsmathematik, VADOR E105-04
\newline TU Wien, Wiedner Hauptstra{\ss }e 8, A-1040 Wien\medskip
\newline\noindent E-mail: \{\texttt{aris.daniilidis, minh.le\}@tuwien.ac.at}
\newline\noindent\texttt{https://www.arisdaniilidis.at/}

\medskip

\noindent Research supported by the Austrian Science Fund grant \textsc{FWF
P-36344N}.\newline\vspace{0.4cm}

\noindent David SALAS

\medskip

\noindent Instituto de Ciencias de la Ingenieria, Universidad de
O'Higgins\newline Av. Libertador Bernardo O'Higgins 611, Rancagua, Chile
\smallskip

\noindent E-mail: \texttt{david.salas@uoh.cl} \newline\noindent
\texttt{http://davidsalasvidela.cl} \medskip

\noindent Research supported by the grants: \smallskip\newline\textsc{CMM
FB210005 BASAL} funds for centers of excellence (\textsc{ANID}-Chile)\newline%
\textsc{FONDECYT 11220586} (Chile)

\end{document}